\documentclass[11pt]{article}
\usepackage{amsfonts}

\usepackage{graphics}

\usepackage{cite}
\usepackage{latexsym}
\usepackage{amsmath}
\usepackage{amssymb}
\usepackage[dvips]{epsfig}
\usepackage{amscd}

\newtheorem{theorem}{Theorem}[section]
\newtheorem{remark}{Remark}[section]

\newtheorem{definition}{Definition}[section]
\newtheorem{lemma}[theorem]{Lemma}

\newcommand{\rt}{\rightarrow}
\newcommand{\vp}{\varphi_{\ve}}
\newcommand{\n}{\rho}
  \newcommand{\nv}{\rho_\ve}
   \newcommand{\uv}{u_\ve}
  
  \newcommand{\ioo}{\int_0^T\int }
  \newcommand{\io}{\int_0^T\int_{\Omega} }
\newcommand{\ti}{\tilde}

\newcommand{\ro}{\rightarrow}

\renewcommand{\div}{ {\rm div }  }

\newcommand{\vf}{v_\ve}
\newcommand{\pa}{\partial}
\renewcommand{\r}{\mathbb{R}}

\renewcommand{\b}{Q_{\ve }}

\newcommand{\bt}{\begin{theorem}}
\newcommand{\bl}{\begin{lemma}}
\newcommand{\el}{\end{lemma}}
\newcommand{\et}{\end{theorem}}
\newcommand{\ga}{\gamma}

\newcommand{\al}{\alpha}
\newcommand{\de}{\delta}
\newcommand{\ve}{\varepsilon}
\newcommand{\la}{\label}

\newcommand{\ol}{\overline}

\newcommand{\bn}{\begin{eqnarray}}
\newcommand{\en}{\end{eqnarray}}
\newcommand{\bnn}{\begin{eqnarray*}}
\newcommand{\enn}{\end{eqnarray*}}

\newcommand{\bnnn}{\begin{eqnarray*}}
\newcommand{\ennn}{\end{eqnarray*}}

\newcommand{\ba}{\begin{aligned}}
\newcommand{\ea}{\end{aligned}}
\newcommand{\be}{\begin{equation}}
\newcommand{\ee}{\end{equation}}
\def\O{{\Omega }}

\def\p{\partial}
\def\norm[#1]#2{\|#2\|_{#1}}

\newcommand{\si}{\sigma}

\def\la{\label}

\def\na{\nabla}

\makeatletter      
\@addtoreset{equation}{section}
\makeatother

\makeatletter      
\@addtoreset{equation}{section}
\makeatother       

  \title{Global Existence of Weak Solutions to the Barotropic Compressible Navier-Stokes Flows  with Degenerate Viscosities}

\date{}

\author{Jing L{\small I}$^{a,b},$ Zhouping  X{\small IN}$^b$\thanks{This research is supported in part by   Zheng Ge Ru Foundation, Hong
Kong RGC Earmarked Research Grants CUHK-4041/11P, CUHK-4048/13P, a Croucher Foundation \& CAS Joint Grant, NSFC/RGC Joint Research Scheme N-CUHK443/14, and a Focus Area Grant at The Chinese University of Hong Kong.
  The research of \textsc{J. Li} was
partially supported  by the National Center for Mathematics and Interdisciplinary Sciences, CAS, and NNSFC  Grant No. 11371348.
 Email:   ajingli@gmail.com (J. Li),
 zpxin@ims.cuhk.edu.hk (Z. Xin).}
 \\
   {\normalsize  a.   Institute of Applied Mathematics, AMSS,}\\
   {\normalsize \&   Hua Loo-Keng Key Laboratory of Mathematics,} \\
  {\normalsize Chinese Academy of Sciences,} \\
   {\normalsize Beijing 100190,
 P. R. China  }\\
  {\normalsize   b.  The Institute of Mathematical Sciences,} \\
   {\normalsize  The Chinese University of Hong Kong, Hong
  Kong}}

\begin{document}
 \maketitle

 \begin{abstract}   This paper concerns  the existence of global  weak solutions to   the barotropic compressible   Navier-Stokes equations with degenerate viscosity coefficients.  We construct suitable approximate system which has smooth solutions satisfying the energy inequality, the BD entropy one, and the Mellet-Vasseur type estimate.  Then,  after adapting  the compactness results due to Bresch-Desjardins (2002, 2003) and Mellet-Vasseur  (2007), we obtain the  global existence of  weak solutions to the  barotropic compressible   Navier-Stokes equations with degenerate  viscosity coefficients      in two or three dimensional periodic domains or whole space for large initial data.  This, in particular, solved an open problem proposed by   Lions (1998).
 \end{abstract}



\section{Introduction and main results}

The barotropic compressible   Navier-Stokes equations, which are the
basic models describing the evolution of a viscous compressible
fluid,
  read as follows \be
\la{ii1}\begin{cases}\rho_t+\div(\rho
u)=0,\\
 (\rho u)_t+\div(\rho u\otimes u)- \div  \mathbb{S} +\nabla P (\rho)=0 ,\end{cases}\ee  where $x\in\O\subset \mathbb{R}^N (N=2,3),
t>0,$ $\n$ is the density, $u=(u_1,\cdots,u_N)$ is the velocity, $\mathbb{S} $ is the viscous stress tensor, and $
P(\rho)=a\rho^\ga(a>0,\ga> 1)$ is the pressure. Without   loss of generality, it is  assumed that $a=1.$
Two major cases will be considered:  either
\be\la{vaa'1} \ba \mathbb{S}\equiv\mathbb{S}_1\triangleq h \na u+g\div u\mathbb{I},\ea\ee or
\be\la{aa'1} \ba \mathbb{S}\equiv\mathbb{S}_2\triangleq h \mathcal{D} u+g\div u\mathbb{I},\ea\ee where $ \mathcal{D}u =\frac{1}{2}(\na u+(\na u)^{\rm tr} ),$  $\mathbb{I}$ is the identical matrix, and $h,g$ satisfy the physical
restrictions \be\la{bd2} h>0,\quad h+Ng\ge 0.\ee

There are many
  studies on the global existence and behavior of solutions
to \eqref{ii1} \eqref{vaa'1} when both  $h$ and $g$ are  constants. The one-dimensional problem
has been studied extensively, see
\cite{ks1,Kan,ho1} and the references therein.
 For the
multi-dimensional case,   the global classical solutions with the density strictly away from vacuum were
first obtained by Matsumura-Nishida \cite{ma2} for initial data
close to a non-vacuum equilibrium. Recently, Huang-Li-Xin \cite{hlx1} obtained the    global classical solutions with the density  containing vacuum provided the initial energy is suitably small. For the weak solutions,
 Hoff \cite{ho3,ho11,ho13}
studied the problem for discontinuous initial data.  When the initial total energy is finite (which implies that the initial density may vanish),  Lions \cite{Lions98}   obtained the
global existence of weak solutions  provided the exponent $\ga$ is suitably large,   which was further relaxed   by Feireisl-Novotny-Petzeltov\'{a}  \cite{fe1} to $\ga>3/2 $  for  three-dimensional case.  

 On the other hand, there are important and interesting phenomena where $h$ and $g$ depend on the density which are degenerate at vacuum. Indeed, as pointed out by  Liu-Xin-Yang in \cite{LiuXinYang98}, in the  derivation of  the
compressible Navier-Stokes equations from the Boltzmann equation
by  the Chapman-Enskog expansions, the viscosity  depends on the
temperature, which is translated into the dependence of the
viscosity on the density for barotropic flows. Moreover, Lions \cite{Lions98} also  proposed various models for shallow water, in particular, he points out  that the global existence of  weak solutions to \eqref{ii1} \eqref{vaa'1} with $h=\n,g=0$ remains open. Recently, 
a friction shallow-water system, with flat bottom topography, which is  derived in \cite{GerbeauPerthame01,mar,brn1,brn2},  can be written in   a two-dimensional space domain $\O$ as \eqref{ii1} \eqref{aa'1} with $h(\n)=g(\n)=\n .$ Indeed, such models appear
naturally and often in geophysical flows
\cite{BreschDesjardins03,BreschDesjardins02,BreschDesjardinsLin05,brn1,brn2}. Therefore, it is of great importance  to study the compressible Navier-Stokes equations  \eqref{ii1} \eqref{vaa'1} and \eqref{ii1} \eqref{aa'1} with density-dependent viscosity.

In the one-dimensional case   with  $h=g=A\n^\al$
for some positive constants $A$ and $\al,$   the well-posedness of either the   initial value problem or the initial  boundary value ones with  fixed or free boundaries has been studied by many authors (see \cite{LiuXinYang98, OkadaMatsusuMakino02,Yangzhao02, Yangzhu02,vong, yangyao01,JiangXinZhang06,llx06} and the references therein).
In higher dimensions, assuming that $h$ is a constant
and $g(\n)=a\n^\beta $ with $a>0 $ and $\beta> 3, $ Vaigant-Kazhikhov \cite{vk95} first proved that for the two-dimensional case \eqref{ii1} \eqref{aa'1} with
slip boundary conditions has a unique global strong and classical
solution. Recently, for the Cauchy problem  and the periodic boundary conditions, Huang-Li \cite{hl1,hl2} and Jiu-Wang-Xin \cite{jwx1,jwx2,jwx3}  relaxed the condition $\beta>3$    to $\beta>4/3.$    For the case $h=h(\n)$ and $g=g(\n),$ in addition to \eqref{bd2}, under the
condition that \be\la{bd1q} g(\n)=h '(\n)\n-h(\n) ,\ee  Bresch-Desjardins \cite{BreschDesjardinsLin05,BreschDesjardins03,
BreschDesjardins03b,BreschDesjardins02} have made important progress. Indeed, for the
periodic boundary conditions and the Cauchy problem, they succeeded in
obtaining a new  entropy inequality (called BD entropy) which can not only be applied
to the vacuum case but also be used to get   the global existence of weak solutions   to \eqref{ii1} \eqref{vaa'1} and  \eqref{ii1} \eqref{aa'1} with some additional drag terms \cite{BreschDesjardinsLin05,BreschDesjardins03,
BreschDesjardins03b}. Later,  by obtaining a new apriori  estimate  on smooth  approximation solutions,  Mellet-Vasseur \cite{MelletVasseur05}    study the stability of barotropic compressible Navier-Stokes equations \eqref{ii1} \eqref{vaa'1} and  \eqref{ii1} \eqref{aa'1} without any additional drag term.   
 However,   the
construction  of the smooth approximation solutions  remains to be carried
out, which does not seem routine in the case of appearance of
vacuum. In fact, only part results for special cases are available.   In particular, for one-dimensional case, Li-Li-Xin \cite{llx06} obtained the global existence of weak solutions to  \eqref{ii1} \eqref{aa'1} with $h(\n)=g(\n)=\n^\al (\al>1/2)$ and proved that for any global entropy weak solution, any vacuum state must vanish within the finite time. Later, when the initial data is   spherically
symmetric,
  Guo-Jiu-Xin \cite{gjx1} obtained the global existence of weak solutions to \eqref{ii1} \eqref{vaa'1} whose Lagrange structure and dynamics   are studied by Guo-Li-Xin \cite{glx}.  Thus, the main aim of this paper is to
 obtain the global existence of weak solutions to  \eqref{ii1} \eqref{vaa'1} and  \eqref{ii1} \eqref{aa'1}    for $\ga>1 $ and for general initial data  by constructing   some suitable  smooth approximation solutions.



For  the sake of simplicity, it is assumed that   for constant $\al>0,$  \be \la{hgvv1}h(\n)=\n^\al,\quad g(\n)=(\al-1)\n^\al.\ee
We  then consider the Cauchy problem,  $\O=\r^N (N=2,3),$  and the case of
bounded domains with periodic boundary conditions,  $\O=\mathbb{T}^N (N=2,3).$   The initial conditions are imposed as \be \la{en1} \n(x,t=0)=\n_0,\quad \n u(x,t=0)=m_0.\ee  We always assume that the initial data $\n_0,m_0$ satisfy that for some  constant  $\eta_0>0,$  \be\la{pini1} \begin{cases} \n_0\ge 0 \mbox{ a.e. in }\O,\,\n_0\not\equiv 0,\,\n_0 \in L^1(\O)\cap L^\ga(\O),\, \na\n_0^{\al-1/2}\in L^2(\O),\\ m_0\in L^{2\ga/(\ga+1)}(\O), \, m_0= 0 \mbox{ a.e. on }\O_0,\,\\
 \n_0^{-1 -\eta_0 } |m_0|^{2+\eta_0}   \in L^1(\O),\end{cases}\ee where we agree that $\n_0^{-1 -\eta_0 } |m_0|^{2+\eta_0} =0$ a.e. on $\O_0,$ the vacuum set of $\n_0 ,$   defined by \be \la{oox1}\O_0\triangleq\{x\in \O\,|\n_0(x)=0\}.\ee 

Before stating the main results, we give the definition of a  weak solution to \eqref{ii1} \eqref{aa'1} \eqref{hgvv1} \eqref{en1}. Similarly, one can define a  weak solution to \eqref{ii1} \eqref{vaa'1} \eqref{hgvv1} \eqref{en1}.
\begin{definition}\la{def1} For $N=2,3,$ let $\O=\mathbb{T}^N $ or $\O=\mathbb{R}^N .$   $(\n,u)$ is said to be a weak solution to \eqref{ii1} \eqref{aa'1} \eqref{hgvv1} \eqref{en1}   if \bnn \begin{cases}0\le \n\in L^\infty(0,T;L^1(\O)\cap L^\ga(\O)),\\ \na\n^{(\ga+\al-1)/2}\in L^2(0,T;(L^2(\O))^N),\\ \na\n^{\al-1/2},\,\,\sqrt{\n}u\in L^\infty(0,T;(L^2(\O))^N),  \\  h(\n) \na u,\,\,h(\n) (\na u)^{\rm tr}\in  L^2(0,T;(W^{-1,1}_{\rm loc}(\O))^{N\times N}),\\ g(\n)  \div u \in  L^2(0,T; W^{-1,1}_{\rm loc}(\O)) ,  \end{cases}\enn
with $ (\n,\sqrt{\n}u)$ satisfying \be\la{fin1} \begin{cases} \n_t+\div (\sqrt{\n}\sqrt{\n}u)=0,\\ \n(x,t=0)=\n_0(x),\end{cases} \mbox{ in }\mathcal{D}', \ee and if the following equality holds for all smooth test function $\phi(x,t)$ with compact support such that $\phi(x,T)=0:$ \be\la{fin2}\ba &\int_\O m_0\cdot \phi(x,0)dx+\io \left(\sqrt{\n}(\sqrt{\n}u)\phi_t+\sqrt{\n}u\otimes \sqrt{\n}u:\na\phi + \n^\ga \div \phi \right)dxdt\\& - \frac12\langle  h(\n) \na u,\na \phi\rangle- \frac12\langle  h(\n) (\na u)^{\rm tr},\na\phi\rangle-  \langle  g(\n) \div u,\div\phi\rangle =0,\ea\ee where \be\nonumber\ba \langle  h(\n) \na u,\na \phi\rangle=&-\io \n^{\al-1/2}\sqrt{\n}u\cdot\Delta \phi dxdt\\&-\frac{2\al}{2\al-1}\io\sqrt{\n}u_j\pa_i\n^{\al-1/2}\pa_i\phi_jdxdt, \ea\ee \be\nonumber\ba \langle  h(\n) (\na u)^{\rm tr},\na\phi\rangle=&-\io \n^{\al-1/2}\sqrt{\n}u\cdot\na\div \phi dxdt\\&-\frac{2\al}{2\al-1}\io\sqrt{\n}u_i\pa_j\n^{\al-1/2}\pa_i\phi_jdxdt, \ea\ee\be\nonumber\ba \langle  g(\n)  \div u ,\div\phi\rangle=&-(\al-1)\io \n^{\al-1/2}\sqrt{\n}u\cdot\na\div \phi dxdt\\&-\frac{2\al(\al-1)}{2\al-1}\io\sqrt{\n}u\cdot\na\n^{\al-1/2}\div\phi dxdt. \ea\ee
\end{definition}

Then the first main result of this paper is as follows:
\begin{theorem}\la{th2}   Let $\Omega=\r^2 $ or $ \mathbb{T}^2 .$  Suppose that  $\al$ and $\ga$ satisfy \be \la{ab1}\al>1/2 ,\quad  \ga>1,\quad\ga\ge 2\al-1. \ee Moreover, assume that the initial data $(\n_0,m_0)$ satisfy \eqref{pini1}. Then  there exists a global weak solution $(\n,u )$ to   the problem \eqref{ii1} \eqref{aa'1} \eqref{hgvv1} \eqref{en1}.

\end{theorem}

 The method of Theorem \ref{th2}   can be applied directly to the system  \eqref{ii1} \eqref{vaa'1}, that is
\begin{theorem}\la{th2'} Let $\Omega=\r^2 $ or $ \mathbb{T}^2 .$  Under the conditions of Theorem \ref{th2},    there exists a global weak solution $(\n,u )$ to   the problem \eqref{ii1} \eqref{vaa'1} \eqref{hgvv1} \eqref{en1}.

\end{theorem}

Theorems \ref{th2} and \ref{th2'} are concerning with the two-dimensional case. As for the three-dimensional case, we have

\begin{theorem}\la{vth2}   Let $\Omega=\r^3 $ or $ \mathbb{T}^3 .$ Suppose that $\al\in [3/4,2) $  and  $\ga\in (1,3)  $ satisfy    \be \ga\in \begin{cases} (1,6\al-3) , &\mbox{ for } \al\in [3/4,1]  ,\\ [2\al-1,3\al-1], &\mbox{ for } \al\in (1,2). \end{cases}  \ee  Assume that the initial data $(\n_0,m_0)$ satisfy \eqref{pini1}.    Moreover,  if  $\al \in(1,2),$ in addition to \eqref{pini1}, we   assume that
\be \la{vin2}\n_0^{-3}|m_0|^4  \in L^1(\O),\ee where we agree that $\n_0^{-3} |m_0|^{4} =0$ a.e. on $\O_0 $    as in \eqref{oox1}. Then there exists a global weak solution $(\n,u )$ to   the problem \eqref{ii1} \eqref{vaa'1} \eqref{hgvv1} \eqref{en1}.
\end{theorem}

\begin{theorem}\la{qvth2}   Let $\Omega=\r^3 $ or $ \mathbb{T}^3 .$ Suppose that $\al=1 $  and  $\ga\in (1,3).  $  Assume that the initial data $(\n_0,m_0)$ satisfy \eqref{pini1}.     Then there exists a global weak solution $(\n,u )$ to   the problem \eqref{ii1} \eqref{aa'1} \eqref{hgvv1} \eqref{en1}.
\end{theorem}

A few remarks are in order:
\begin{remark}    If $\al=1$ and $\gamma=2,$ Theorems \ref{th2}--\ref{qvth2} give a positive answer to   the open problem proposed by Lions   \cite[Section 8.4]{Lions98}: ``In the first case (ie (8.70)--(8.71)), the Cauchy problem is completely open for the models involving  (8.73)", where  (8.70)--(8.71) (8.73) is  corresponding to \eqref{ii1} \eqref{vaa'1} \eqref{hgvv1}  with $\al=1$ and $\ga=2.$  \end{remark}

\begin{remark}   For three-dimensional case, it should be noted that  Theorem  \ref{vth2} here is valid for all $\ga\in (1,3)$ provided $h=\n$ and $g=0.$ Therefore, for  $h=\n$ and $g=0,$  our Theorems \ref{vth2} and \ref{qvth2}  establish the existence of global weak solutions to   \eqref{ii1} \eqref{vaa'1} and  \eqref{ii1} \eqref{aa'1}   with   $\ga\in (1,3) $ for general initial data, which is in sharp contrast to   the case that $h$ and $g$ are  both constants, where the condition $\ga>3/2 $  is essential in the
analysis of  Lions \cite{Lions98} and Feireisl-Novotny-Petzeltov\'{a} \cite{fe1}.  In fact, for $h$ and $g$  being  both constants and $\ga\in (1,3/2],$  it remains completely  open to obtain the global existence of weak solutions to \eqref{ii1} \eqref{vaa'1} for general initial data  except for the spherically symmetric case \cite{jz1}.
\end{remark}

\begin{remark}    After some routine modifications, for the system  \eqref{ii1} \eqref{vaa'1}, our method can be applied directly to the case that  $h$ and $g$ satisfy  \eqref{bd1q} and the conditions listed in \cite{MelletVasseur05} together with some additional    constraints. However,
 for the shallow-water system  with $h=g=\n $ (\cite{GerbeauPerthame01,mar,brn1,brn2}), the global existence of weak solutions to \eqref{ii1} \eqref{vaa'1} or \eqref{ii1} \eqref{aa'1} for general initial data  remains open since \eqref{bd1q} fails for this case. \end{remark}

\begin{remark}  For the system  \eqref{ii1} \eqref{aa'1}  and  three-dimensional case, our construction depends on the condition that $h=\rho$ and $g=0$,  and cannot be applied directly to the general case that $h$ and $g$ satisfy the conditions listed in Theorem \ref{vth2}. This will be left for future. \end{remark}

\begin{remark}  Around the same time when this paper is finished, there are announcements of some existence results on the  problem \eqref{ii1} \eqref{aa'1} \eqref{hgvv1} \eqref{en1} with $ \al=1 $ and $\O=\mathbb{T}^N (N=2,3) $ by   Vasseur-Yu \cite{vy} with a different approach. However, we have difficulties to understand some of their key a priori assumptions near vacuum in their arguments.
\end{remark}

We now make some comments on the analysis of this paper.
Since the compactness arguments are similar  to those of Bresch-Desjardins \cite{BreschDesjardins03,BreschDesjardins03b, BreschDesjardinsLin05} and Mellet-Vasseur \cite{MelletVasseur05}, the main point of this paper is to construct smooth approximate  solutions, whose densities are  bounded  from above and strictly bounded away from vacuum provided the smooth initial ones are,  satisfying the energy estimate, the BD entropy inequality, and the  Mellet-Vasseur type estimate.  To this end, we first deal with the periodic case and consider the following approximate  system
     \be\la{ba1} \begin{cases}   \n_t+\div(\n u)  =   \ve \n^{1/2}\div(\n^{- 1/2}h_\ve'(\n)\na \n),
     \\   \n u_t+ \n u\cdot \na u   -\div(h_\ve(\n)   \mathcal{D} u )  -\na(g_\ve(\n)\div u)  +\na P\\ =\sqrt{\ve}\div(h_\ve(\n)    \na u) +\sqrt{\ve}\na(g_\ve(\n)\div u)-e^{-\ve^{-3}}  (\n^{ \ve^{-2}}+ \n^{-\ve^{-2}})u ,  \end{cases}\ee where
\be \la{ini9} h_\ve(\n)=\n^\al+\ve^{1/3} ( \n^{7/8} + \n^{\ti\ga }), \quad g_\ve(\n)=  \n h_\ve'(\n)-h_\ve(\n),\ee
with \be\la{ma1}  0<\ve\le \ve_0\triangleq \min\{(2\al-1) (16(\al+\ga))^{-10},\eta_0\},\quad \ti\ga\triangleq \ga+1/6.\ee
Here, we propose to approximate \eqref{ii1}$_1$ by \eqref{ba1}$_1$ which is a parabolic equation for any fixed $\ve>0$ and hence has   smooth effects  on the density provided the smooth initial density is strictly away from vacuum. The specific choices  of the higher order regularization in \eqref{ba1} have several key advantages. First, it can be shown that the smooth solutions to the  new system \eqref{ba1}--\eqref{ma1} satisfy the energy   and the  Mellet-Vasseur type estimates. Moreover,   after some careful calculations, we find that the most difficult term   induced by $\ve \n^{1/2}\div(\n^{- 1/2}h_\ve '(\n)\na \n)$ has the right sign (see \eqref{bb8}) which implies that the solutions  to our approximate system also satisfy the BD entropy inequality. In fact, this is one of the key observations of this paper.  Next, in order to obtain the lower and upper bounds of the density, in addition to the   estimate on  $L^\infty(0,T;L^{N+\de})$-norm of $\n^{1/(N+\de)} u$ which can be obtained for the system  \eqref{ii1} \eqref{aa'1}   in two-dimensional case (see \eqref{e1}) and for   \eqref{ii1} \eqref{vaa'1} in both two-dimensional (see \eqref{e1}) and three-dimensional  cases (see \eqref{k1}),  one still needs some additional estimates on the $L^\infty(0,T;L^p)$-norm (for suitably large   $p$) of $\n$ and $\n^{-1}$ which can be achieved by adding a damping term $-e^{-\ve^{-3}}  (\n^{ \ve^{-2}}+ \n^{-\ve^{-2}})u $ on the righthand side of \eqref{ii1}$_2$ (see \eqref{ba1}$_2$). However, for $\ve\to 0^+,$ this term will bring new difficulties   which can be overcome by adding  $\ve^{1/3}(\n^{7/8}+\n^{\ga+1/6})$ to $h(\n)$ (see   \eqref{ini9}). This idea is motivated by our previous study on the one-dimensional problem \cite{llx06}.
With all these estimates at hand, we can use   a De Giorgi-type procedure to bound the density from above and below, in particular, the density is strictly away from vacuum provided the initial one is (see \eqref{c1}). In fact, this is another key issue of this paper. Once we obtained \eqref{c1}, we can use the $L^p$-theory for parabolic system to get the estimates on the $L^p(0,T;L^p)$-norm of $(\n,u),(\n_t,u_t),$ and $(\na^2\n,\na^2u)$ (see \eqref{bc4}). This in turn implies that  the approximate system \eqref{ba1}--\eqref{ma1}  has a global strong solution with smooth initial data. Next, after adapting   the compactness results due to Bresch-Desjardins \cite{BreschDesjardins03,BreschDesjardins03b,BreschDesjardinsLin05} and Mellet-Vasseur \cite{MelletVasseur05}, we can  obtain the global existence of the weak solutions to either \eqref{ii1} \eqref{aa'1} \eqref{hgvv1} for two-dimensional periodic case or \eqref{ii1} \eqref{vaa'1} \eqref{hgvv1} for two-dimensional and  three-dimensional  periodic  cases.  Finally, to prove Theorems \ref{th2} and \ref{vth2}  where $\O=\r^N (N=2,3), $   some     extra care should be taken due to the unboundedness of the domain.    In this case, we consider the system \eqref{ba1}--\eqref{ma1} in $\b=(-\ve^{-\si_0},\ve^{-\si_0})^2  $  and \eqref{vba1} \eqref{ini9} \eqref{ma1} in $\b=(-\ve^{-\si_0},\ve^{-\si_0})^N (N=2,3)$ and impose the Neumann boundary condition on $\n$ and Navier-slip conditions on   $u$ (see \eqref{vbb35} and \eqref{vbb35v}). Then we can adapt the preceding proofs in the case $\O=\mathbb{T}^N  (N=2,3)$ to  $\O=\r^N (N=2,3).$

This paper is organized as follows. Since  the proof of Theorem \ref{th2'} is similar as that of Theorem \ref{th2} after some routine modifications,  we will only prove Theorems \ref{th2},   \ref{vth2}, and \ref{qvth2}.  In the next section,
we   work   on the problem \eqref{ii1} \eqref{aa'1} \eqref{hgvv1}  \eqref{en1} in the two-dimensional periodic case, $\O=\mathbb{T}^2,$ then in the Section 3, we adapt the previous procedure to the problem \eqref{ii1} \eqref{vaa'1} \eqref{hgvv1} \eqref{en1}   in the three-dimensional periodic case, $\O=\mathbb{T}^3.$ Next, in the section 4, we will construct a new approximate system which can be applied to obtain the global weak solutions to the problem \eqref{ii1} \eqref{aa'1} \eqref{hgvv1} \eqref{en1}   in the three-dimensional periodic case,   and in the end (Section 5) we shall explain how to modify the preceding proofs in the cases where $\O=\r^N (N=2,3) .$

\section{Proof  of Theorem  \ref{th2}: $\O=\mathbb{T}^2$}

In this section, we study the 2-dimensional periodic case which is the simplest one, yet the most important case since most of ideas to be developed here can be modified to deal with other cases.

\subsection{A priori estimates }
For $\ve$ as in \eqref{ma1}, let smooth  functions $ \n_{0\ve}>0$ and $u_{0\ve} $ satisfy
\be \la{pd9}\ba&\|\n_{0\ve}\|_{L^1\cap L^{\ga}(\O) } +\|\na  \n_{0\ve}^{\al-1/2 } \|_{L^2(\O)}+\ve^{1/3}\|\na  \n_{0\ve}^{3/8 } \|_{L^2(\O)}+\ve^{1/3}\|\na  \n_{0\ve}^{\ga-1/3 } \|_{L^2(\O)}\\&+\ve^{13/3}e^{-\ve^{-3}}   \|  \n_{0\ve}^{ \ve^{-2}+\ti\ga -1} \|_{L^1(\O)} +\ve^{13/3}e^{-\ve^{-3}}   \|  \n_{0\ve}^{-\ve^{-2}-1/8} \|_{L^1(\O)}\\&+\ve^{4}e^{-\ve^{-3}}   \|  \n_{0\ve}^{ \ve^{-2}+\al -1} \|_{L^1(\O)} +\ve^4e^{-\ve^{-3}}   \|  \n_{0\ve}^{-\ve^{-2}+\al -1} \|_{L^1(\O)}\le C, \ea \ee and \be  \la{pd22}\int_\O \n_{0\ve}   |u_{0\ve}|^{2+\eta_0}dx\le C,\ee for some constant $C$ independent of $\ve.$  We extend  $\n_{0\ve} $ and $u_{0\ve}$ $\O$-periodically to $\r^2 $ and  consider  the system \eqref{ba1}--\eqref{ma1} with initial data:
\be \la{bb35} (\n,u)(x,0)= (\n_{0\ve},u_{0\ve})  .\ee
 Let $T>0$ be a fixed time and $(\rho ,u)$   be
a smooth solution to \eqref{ba1}--\eqref{ma1} \eqref{bb35}  on
$\O \times (0,T]. $

Then, we will establish some necessary a priori bounds
for  $(\rho,u)$. The first one is the energy-type inequality.

\begin{lemma} \la{lem10} There exists some generic constant $C$ independent of $\ve$ and $T$ such that
\be\la{bb30}\ba& \sup_{0\le t\le T}\int(\n |u|^2+\n+  \n^\ga  )dx    +\ve\int_0^T\int    \n^{-1} h_\ve'(\n) |\na \n|^2(1+|u|^2) dxdt \\&+\int_0^T\int  h_\ve(\n)  |\mathcal{D} u|^2 dxdt+ e^{-\ve^{-3}} \int_0^T    \int(\n^{ \ve^{-2}}+\n^{-\ve^{-2}})|u|^2dx  dt  \le C,\ea\ee where and throughout  this section, for any $f$,
\bnn \int fdx\triangleq\int_\O fdx.\enn
\end{lemma}

{\it Proof.} First, integrating $\eqref{ba1}_1$ over $\O\times (0,T)$ together with \eqref{bb35} gives \be \la{bc6}\sup_{0\le t\le T}\int\n dx+\ve\int_0^T\int \n^{-1}h'_\ve(\n)|\nabla \n|^2dxdt\le C.\ee

Next, multiplying $\eqref{ba1}_2$ by $u,$  integrating by parts, and using $\eqref{ba1}_1$ yield
\be\la{bb61}\ba& \frac12 (\int\n |u|^2dx)_t   +\int (h_\ve(\n) (|\mathcal{D}u|^2+\sqrt{\ve}|\na u|^2)+(1+\sqrt{\ve}) g_\ve(\n) (\div u)^2)dx \\&\quad +e^{-\ve^{-3}}   \int\left(\n^{ \ve^{-2}}+\n^{-\ve^{-2}}\right)|u|^2dx  +\int   u \cdot \na\n^\ga  dx\\&=\frac{\ve}{2}\int   \n^{1/2}\div(\n^{- 1/2}h_\ve'(\n)\na \n) |u|^2dx \\&=-\frac{\ve}{4} \int \n^{-1}h_\ve'(\n)  |\na \n |^2|u|^2dx-  \ve \int h_\ve'(\n)   \na \n \cdot\na u\cdot udx\\&\le-\frac{\ve}{8} \int \n^{-1}h_\ve'(\n)  |\na \n |^2|u|^2dx+2\ve \int\n h_\ve'(\n)   |\na u|^2dx  \\&\le-\frac{\ve}{8} \int \n^{-1}h_\ve'(\n)  |\na \n |^2|u|^2dx+\frac{\sqrt{\ve}}{2} \int  h_\ve (\n)   |\na u|^2dx    .\ea\ee

Then, to estimate the last term on the left hand side  of \eqref{bb61}, after integration by parts and using $\eqref{ba1}_1$, one obtains  that for $q\not=1,$
\be\la{bini5}\ba    \int  u\cdot\na  \n^{q }  dx    &= -\frac{q}{q -1}\int\n^{q -1}\div(\n u) dx\\ &= -\frac{q}{q -1}\int\n^{q -1}(-\n_t + \ve \n^{1/2}\div(\n^{- 1/2}h_\ve'(\n)\na \n)   ) dx \\ &= \frac{1}{   q-1 }(\int\n^{q }dx)_t+\frac{q(2q -1)\ve}{2(q -1)}\int\n^{q -2} h_\ve'(\n)|\na \n|^2dx .\ea\ee

 Finally, for $v\in \r^N (N=2,3),$  we have  $$ (\div v)^2  \le N  |\mathcal{D} v|^2 \le N|\na v|^2 , $$  which together with \eqref{ini9}  implies that for $N=2,3,$ \be \la{bb63} \begin{cases}4(h_\ve(\n) |\mathcal{D}v  |^2 +   g_\ve(\n) (\div v)^2 )\ge
  \min\{N\al-(N-1),1\}  h_\ve(\n)  |\mathcal{D} v|^2,\\ 4(h_\ve(\n) |\na v|^2 +   g_\ve(\n) (\div v)^2 )\ge
  \min\{N\al-(N-1),1\}  h_\ve(\n)  |\na v|^2.\end{cases}\ee    Since $\ve\le \ve_0,$ the combination of  \eqref{bc6}--\eqref{bb63}, \eqref{pd9}, with \eqref{pd22}  yields \eqref{bb30}, which completes the  proof of Lemma \ref{lem10}.

Now we are in a position to derive the following  entropy estimate which in particular yields the uniform BD  one due to Bresch-Desjardins \cite{BreschDesjardinsLin05,BreschDesjardins03,
BreschDesjardins03b,BreschDesjardins02}.

\begin{lemma}
\la{lem11} There  exists some generic constant $C$ independent of $\ve$ and $T$ such that
\be\la{bb5'}\ba& \sup_{0\le t\le T}\int  \left( \n^{-1}(h_\ve'(\n))^2|\na   \n  |^2  +\ve^{13/3}e^{-\ve^{-3}}    (\n^{ \ve^{-2}+\ti\ga-1 }+ \n^{-\ve^{-2}-1/8})   \right)dx \\& + \int_0^T \int h_\ve(\n)  |\na u|^2dxdt  + \int_0^T \int \n^{\ga-3}h_\ve(\n)  |\na \n|^2dxdt  \le C.\ea\ee

\end{lemma}

{\it Proof. } First, set
\be  \la{lalm1} G  \triangleq  \ve  \n^{1/2}\div(\n^{- 1/2}h_\ve'(\n)\na \n)
\ee and
\be \la{bini4} \vp'(\n)\triangleq \n^{-1}h_\ve'(\n)\ge 0
   .\ee
Multiplying \eqref{ba1}$_1$ by $\vp '(\n)$ leads to
\be\nonumber\la{bb1} (\vp( \n))_t+u\cdot\na \vp( \n)  + \n\vp'(\n) \div u= \vp'(\n)  G,\ee which gives
 \be\la{bb2} (\na\vp( \n) )_t+u\cdot\na\na\vp( \n)  +\na u\cdot\na \vp( \n) + \na(\n\vp'( \n) \div u)= \na(\vp'( \n) G).\ee
Thus, multiplying \eqref{bb2} by $\n\na\vp (\n)$ and integration by parts show that
 \be\la{bb3} \ba&\frac{1}{2}(\int\n|\na\vp( \n)|^2dx)_t  +\int \na h_\ve(\n)\cdot \na u\cdot\na\vp( \n) dx \\& + \int\na h_\ve(\n)\cdot \na(\n\vp'(\n)\div u)dx\\& +\int\vp'(\n) G\left( \Delta h_\ve(\n)-\frac12\vp'(\n)|\na  \n|^2  \right)dx=0.\ea\ee

Next, multiplying \eqref{ba1}$_2$ by $\na \vp (\n)$  leads to
\be\la{bb4}\ba& \int  u_t\cdot \na h_\ve( \n) dx + \int u\cdot \na u\cdot \na h_\ve(\n) dx-(1+\sqrt{\ve})  \int h_\ve(\n)\na\div u\cdot\na\vp(\n)dx\\&-(1+\sqrt{\ve})\int \na h_\ve( \n)\cdot \na u\cdot \na \vp( \n) dx+(1+\sqrt{\ve})\int g_\ve(\n)\div u\Delta\vp(\n)dx  \\&+ \int P'(\n)\vp'(\n)|\na \n|^2dx+e^{-\ve^{-3}}\int ( \n^{ \ve^{-2}}+\n^{-\ve^{-2}}) u\cdot\na \vp(\n) dx=0,\ea\ee where  the following simple fact  has been used:
\be \nonumber\la{bb64} \ba -\int h_\ve(\n)\Delta u\cdot \na\vp( \n )dx=-\int h_\ve(\n)\na\div u\cdot \na\vp( \n )dx.\ea\ee

Since \eqref{ba1}$_1$ implies  $$ (h_\ve(\n))_t+\div(h_\ve(\n)u)+(\n h_\ve'(\n)-h_\ve(\n))\div u=h_\ve'(\n)G,$$
the first term on the left hand side of \eqref{bb4} is handled as
\be\la{bb10}\ba  \int  u_t\cdot \na h_\ve(\n) dx  &=(\int  u\cdot \na h_\ve( \n) dx )_t-\int  u\cdot \na h_\ve(\n)_t dx 
 \\&=(\int  u\cdot \na h_\ve( \n) dx )_t-\int u\cdot\na  u\cdot \na h_\ve( \n )  dx \\&\quad -2\int h_\ve( \n )\mathcal{D}u:\na u   dx+ \int h_\ve( \n )|\na u|^2   dx\\&\quad -\int(\n h_\ve'(\n)- h_\ve(\n))(\div u)^2 dx +\int \div u h_\ve'(\n)G dx , \ea\ee
where in the second equality one has used
\bnn \ba & \int    u\cdot \na \div ( h_\ve(\n)u) dx \\&=  -\int \pa_i  u\cdot \na(h_\ve( \n )u_i) dx  \\&=   -\int u\cdot\na  u\cdot \na h_\ve( \n )  dx -\int h_\ve( \n )\pa_i  u\cdot \na u_i  dx \\&=  -\int u\cdot\na  u\cdot \na h_\ve( \n )  dx -2\int h_\ve( \n )\mathcal{D}u:\na u   dx+ \int h_\ve( \n )|\na u|^2   dx . \ea\enn

Now, multiplying \eqref{bb3}  by $ 1+\sqrt{\ve} $ and adding the resulting equality to \eqref{bb4}, one can obtain after using \eqref{bb10}  that
\be\la{bb7}\ba& \frac{1+\sqrt{\ve} }{2}(\int\n|\na \vp( \n) |^2dx)_t + (\int \n  u\cdot \na \vp(\n )dx )_t+ \int h_\ve(\n)|\na u|^2dx \\&  \quad +  \int P'(\n)\vp'(\n)|\na \n|^2dx + {e^{-\ve^{-3}}} \int  ( \n^{ \ve^{-2}} + \n^{-\ve^{-2}} ) u\cdot\na \vp( \n) dx\\& \quad+(1+\sqrt{\ve}) \int\vp'(\n) G\left( \Delta h_\ve(\n)-\frac12\vp'(\n)|\na  \n|^2 +\frac{1}{1+\sqrt{\ve} }\n\div u \right)dx\\&=2\int h_\ve( \n)\mathcal{D}u: \na  u   dx + \int(\n h_\ve'(\n)- h_\ve(\n))(\div u)^2 dx\\&\le \frac{1}{2 }\int h_\ve( \n)| \na  u|^2   dx +C\int  h_\ve (\n) |\mathcal{D} u|^2 dx,\ea\ee
where in the first equality one has used the following simple calculations:
\be \nonumber\la{bb65} \ba  & \int\na h_\ve(\n)\cdot \na(\n\vp'(\n)\div u)dx -\int h_\ve(\n)\na\div u\cdot\na\vp(\n)dx\\&\quad+ \int g_\ve(\n)\div u\Delta\vp(\n)dx\\&  =\int\na h_\ve(\n)\cdot \na(\n\vp'(\n))\div u dx+\int \n\vp'(\n)\na h_\ve(\n)\cdot \na \div u dx\\&\quad- \int h_\ve(\n)\na\div u\cdot\na\vp( \n )dx-  \int g_\ve(\n) \na \div u\cdot\na\vp(\n) dx  \\&\quad- \int  \div u \na g_\ve(\n)\cdot\na\vp(\n) dx\\&  =\int\left(\na h_\ve(\n)\cdot \na(\n\vp'(\n))- \na g_\ve(\n)\cdot\na\vp(\n)\right)\div u dx \\&\quad+\int\left( \n\vp'(\n)\na h_\ve(\n) - h_\ve(\n) \na\vp( \n )-   g_\ve(\n) \na\vp(\n)\right)\cdot\na \div u dx=0 \ea\ee
due to \eqref{ini9} and \eqref{bini4}.

Since \eqref{bini4}  and \eqref{lalm1} imply
\bnn  \Delta h_\ve( \n)-\frac12\vp'(\n)|\na  \n|^2    =    \ve^{-1}G,\enn the last term on the left hand side of \eqref{bb7} satisfies
\be\la{bb8}\ba &(1+\sqrt{\ve} )\int\vp'(\n) G\left( \Delta h_\ve(\n)-\frac12\vp'(\n)|\na  \n|^2 +\frac{1}{1+\sqrt{\ve} }\n\div u \right)dx \\&\ge \frac{1}{2\ve}\int  \vp'(\n) G^2dx- \frac{\ve}{2}\int\n^2 \vp'(\n) \left(\div u\right)^2dx \\&\ge \frac{1}{2\ve}\int \vp'(\n) G^2dx-C{\ve}  \int  h_\ve (\n) |\mathcal{D} u|^2dx.\ea\ee

Finally, it follows from \eqref{bini4} and \eqref{ini9} that
\bnn\ba  & \int  (\n^{ \ve^{-2}}+\n^{-\ve^{-2}}) u\cdot\na \vp( \n) dx \\& =  \int  u\cdot\na \left(\frac{\al \n^{ \ve^{-2}+\al-1}}{ \ve^{-2}+\al-1}  +\frac{7\ve^{7/3}\n^{\ve^{-2}-1/8}}{ 8-\ve^2} +\frac{\ti\ga\ve^{7/3} \n^{ \ve^{-2}+\ti\ga -1}}{ 1+(\ti\ga -1)\ve^2}  \right) dx\\&\quad+\int  u\cdot\na \left(\frac{\al \n^{-\ve^{-2}+\al-1}}{-\ve^{-2}+\al-1}  -\frac{7\ve^{7/3}\n^{-\ve^{-2}-1/8}}{ 8+\ve^2} -\frac{\ti\ga\ve^{7/3} \n^{-\ve^{-2}+\ti\ga-1}}{1-(\ti\ga -1)\ve^2}   \right) dx ,\ea\enn which, together with  \eqref{bini5},  \eqref{bb7}, \eqref{bb8},     \eqref{bb30},  and   \eqref{pd9},  yields \eqref{bb5'}. The proof of Lemma \ref{lem11} is finished.

With Lemmas  \ref{lem10} and  \ref{lem11} at hand, we can prove the following   Mellet-Vasseur type estimate (\cite{MelletVasseur05}).
\begin{lemma} \la{lem01} Assume that $\ga>1$ satisfies $\ga\ge (1+\al)/2 $ in addition.     Then  there exists some generic constant $C$    depending on $T$ but independent of $\ve$ such that \be \la{bb'01} \sup_{0\le t\le T}\int\n (e+|u|^2)\ln (e+|u|^2)dx\le C.\ee\end{lemma}

{\it Proof.} First, multiplying \eqref{ba1}$_2$ by
$(1+\ln (e+|u|^2))u $ and integrating lead to
\be \la{bb67} \ba &\frac{1}{2} \frac{d}{dt}\int\n (e+|u|^2)\ln (e+|u|^2)dx-\frac{1}{2} \int (e+|u|^2)\ln (e+|u|^2)G dx\\&\quad+\int (1+\ln (e+|u|^2)) (h_\ve(\n) (|\mathcal{D} u|^2+\sqrt{\ve} |\na u|^2) + (1+\sqrt{\ve}) g_\ve(\n)  (\div u)^2)dx  \\&\le C\int h_\ve(\n) |\na u|^2dx-\int(1+\ln (e+|u|^2))u\cdot\na \n^\ga dx\\&\le C\int h_\ve(\n) |\na u|^2dx+C \int  \ln^2 (e+|u|^2)  \n^{2\ga-\al}  dx  ,\ea\ee
where in the last inequality one has used the following estimate
\be \la{bb1x} \ba &\left|\int(1+\ln (e+|u|^2))u\cdot\na \n^\ga dx\right|\\ &\le \int (1+\ln (e+|u|^2))|\div u|\n^\ga dx +\left|\int \frac{2u_iu_k}{e+|u|^2}\pa_iu_k  \n^\ga dx\right| \\ &\le  C \int  \ln^2 (e+|u|^2)  \n^{2\ga-\al}  dx +C\int h_\ve(\n)|\na u|^2 dx .\ea\ee

Then, integration by parts gives
\be  \la{bb2x} \ba & -\frac{1}{2} \int (e+|u|^2)\ln (e+|u|^2)G dx\\& = \frac{\ve}{4} \int \n^{-1}h_\ve'(\n)  |\na \n |^2(e+|u|^2)\ln (e+|u|^2)dx+ \ve \int h_\ve'(\n)    \na \n \cdot\na u\cdot udx\\&\quad+  \ve \int h_\ve'(\n)  \ln (e+|u|^2) \na \n \cdot\na u\cdot udx\\&\ge \frac{\ve}{8}  \int \n^{-1}h_\ve'(\n)  |\na \n |^2(e+|u|^2)\ln (e+|u|^2)dx-\ve \int\n h_\ve'(\n)   |\na u|^2dx\\& \quad- {\ve}  \int \n^{-1}h_\ve'(\n)  |\na \n |^2|u|^2dx-2\ve (\al+2)\int  h_\ve (\n)  \ln (e+|u|^2)  |\na u|^2dx . \ea\ee
It follows from this, \eqref{bb67},   \eqref{bb30}, \eqref{bb5'},     \eqref{bb63}, and  \eqref{pd22}  that
\be \la{com1}\sup_{0\le t\le T}\int\n (e+|u|^2)\ln (e+|u|^2)dx\le C +\int_0^T \int  \ln^2 (e+|u|^2)  \n^{2\ga-\al}  dx dt.\ee

Finally, since $\ga\ge(\al+1)/2,$  it holds that
 \be \la{bb70} \ba  \int  \ln^2 (e+|u|^2)  \n^{2\ga-\al}  dx  &\le  C \int   (\n+\n^{2\ga})(1+|u| ) dx  \\&\le C+ C \int\n |u|^2dx+C \int(\n +\n^{4\ga-1})dx \\&\le C, \ea\ee where in the last inequality, one has used \eqref{bb30}, \eqref{bb5'}, and the following Sobolev inequality that for any $p>1,$ there exists some constant $C $ depending only on $\al$ and $p$ such that \be \la{bb6'} \|\n\|_{L^p(\O)} \le C \|\n \|_{L^1(\O)}+C \|\na \n^{\al-1/2}\|_{L^2(\O)}^{2/(2\al-1)}  .\ee

Putting \eqref{bb70}  into \eqref{com1} yields \eqref{bb'01}. The proof of Lemma \ref{lem01} is completed.

  Next, we will use  a De Giorgi-type procedure to obtain the following estimates on the lower and upper bounds of the density which are the key  to obtain the global existence of strong solutions to the problem \eqref{ba1}--\eqref{ma1}  \eqref{bb35}.

\begin{lemma}
\la{lem13}There exists some positive constant $C $ depending on $\ve$ and $T$ such that for all $(x,t)\in \O\times (0,T)$ \be \la{c1} C^{-1}\le \n(x,t)\le C . \ee \end{lemma}
{\it Proof.} First, multiplying \eqref{ba1}$_2$ by
$ |u|^\ve u $ and integrating in space give
\be \la{b'67} \ba &\frac{1}{2+\ve} \frac{d}{dt}\int\n  |u|^{2+\ve} dx - \frac{1}{2+\ve} \int  |u|^{2+\ve} G dx  \\&\quad+\int |u|^{\ve}\left(h_\ve(\n) (|\mathcal{D}u|^2+\sqrt{\ve} |\na u|^2) + (1+\sqrt{\ve}) g_\ve(\n)  (\div u)^2\right)dx \\&\quad+\frac{\ve(1+2\sqrt{\ve} )}{2}\int h_\ve(\n)|u|^\ve |\na |u||^2dx +e^{-\ve^{-3}}\int (\n^{\ve^{-2}}+\n^{-\ve^{-2}})|u|^{2+\ve}dx \\&=   -\frac{1}{2}\int h_\ve(\n)u\cdot\na u\cdot\na |u|^\ve dx-(1+\sqrt{\ve})\int g_\ve(\n)  \div u u\cdot\na |u|^\ve dx\\&\quad  -\int |u|^\ve u\cdot\na \n^\ga dx\\&\le 4(\al+\ga)\ve\int h_\ve(\n)|u|^\ve |\na u|^2dx  +C\int h_\ve(\n)|\na u|^2 dx \\& \quad+C \int   \left(\n^{\ve^{-2}}+\n^{-\ve^{-2}}\right)  |u|^{2 }  dx+C,\ea\ee  where in the last inequality one has used
 the following simple fact  that  \be \la{bn1} \sup_{0\le t\le T}\int \left(\n^{\ve^{-2}}+\n^{-\ve^{-2}}\right)dx \le C ,\ee due to \eqref{bb5'}.
Integration by parts yields that
\be \la{bc72} \ba & -\frac{1}{2+\ve} \int  |u|^{2+\ve} G dx\\& = \frac{\ve}{2(2+\ve)} \int \n^{-1}h_\ve'(\n)  |\na \n |^2 |u|^{2+\ve}dx+ \ve  \int h_\ve'(\n)   |u|^{1+\ve}  \na \n \cdot\na |u| dx \\&\ge \frac{\ve }{8(2+\ve)}  \int \n^{-1}h_\ve'(\n)  |\na \n |^2 |u|^{2+\ve} dx-2(2+\ve)\ve \int\n h_\ve'(\n)  |u|^\ve |\na u|^2dx   .\ea\ee
It follows from \eqref{b'67},  \eqref{bc72},    \eqref{bb63}, \eqref{bb30},  \eqref{bb5'}, and \eqref{ma1} that \be\la{e1}\sup_{0\le t\le T}\int\n  |u|^{2+\ve}dx+\sqrt{\ve}\int_0^T\int h_\ve(\n)|u|^\ve |\na u|^2dxdt\le C.\ee

Next,  since $v\triangleq \n^{1/2} $ satisfies
\be\la{bb22}2 v_t - 2\ve \div(h_\ve'(v^2)\na v )+\div(u  v)+ u\cdot\na v=0,  \ee
multiplying \eqref{bb22} by $(v-k)_+$ with $k\ge \|v(\cdot,0)\|_{L^\infty(\O)}  = \| \n_0\|_{L^\infty(\O)}^{1/2} $ and integrating by parts yield
\be\la{bb20}\ba &\frac{d}{dt}\int (v-k)_+^2dx+ 2\al \ve \int v^{2\al-2}|\na(v-k)_+ |^2dx   \\&\le C\int_{A_k(t)}v^{4-2\al}|u|^2 dx+ \al \ve \int v^{2\al-2}|\na(v-k)_+ |^2dx,\ea\ee where  $A_k(t)\triangleq\{x\in\O|v(x,t)>k\} .$
 It thus follows from \eqref{e1} and H\"older's inequality  that
\be\la{bb'20}\ba &\int_{A_k(t)}v^{4-2\al}|u|^2 dx \\ &\le C\left(\int_{A_k(t)}v^2|u|^{2+\ve}dx\right)^{2/(2+\ve)} \left(\int_{A_k(t)}v^{ (4+4\ve-2(2+\ve)\al)/\ve}dx\right)^{\ve/(2+\ve)}\\ &\le C \left(\int_{A_k(t)}(\n^{4(\al+1)\ve^{-1}} +\n^{-4(\al+1)\ve^{-1}})dx\right)^{\ve/(2+\ve)} \\ &\le C \left(\int_{A_k(t)}(\n^{ \ve^{-2}}+\n^{- \ve^{-2}})dx\right)^{\ve(4-\ve)/(6(2+\ve))} |A_k(t)|^{\ve/6}\\ &\le C|A_k(t)|^{\ve/6 },\ea\ee where \eqref{bn1}   has been used in the last inequality. Putting \eqref{bb'20} into \eqref{bb20} leads to
\be\la{bb21}\ba & I_k'(t)+   \al \ve \int \n^{ \al-1}|\na(v-k)_+ |^2dx \le C \nu_k^{\ve/6 } ,\ea\ee
 where
$$ I_k(t)\triangleq\int (v-k)_+^2(x,t)dx,\quad \nu_k\triangleq\sup\limits_{0\le t\le T}|A_k(t)|.$$
   Since $I_k(0)=0,$ without loss of generality, we can assume that there exists some $\si>0$  such that $$I_k(\si)=\sup_{0\le t\le T}I_k(t).$$ It follows from \eqref{bb21} that \bnn\la{bb11}   I_k(\si)+\int \n^{ \al-1}|\na(v-k)_+ |^2(x,\si)dx \le C\nu_k^{\ve/6},\enn which, together with
H\"older's inequality and \eqref{bn1}, gives \be\la{bb12}\ba &I_k(\si)+\|\na  (v-k)_+(\cdot,\si)\|_{L^{24/(12+\ve)}(\O)}^2 \\ &\le C\nu_k^{\ve/6 }+\int \n^{ \al-1}|\na(v-k)_+ |^2(x,\si)dx \left(\int  \n^{12(1-\al)/\ve}(x,\si)dx\right)^{ \ve/12}\\ &\le C\nu_k^{\ve/6 }.\ea\ee

Then, for any $ h>k\ge \|v(\cdot,0)\|_{L^\infty(\O)} ,  $ direct computations yield
\bnn\ba &|A_h(t)|(h-k)^2\\&\le \|(v-k)_+(\cdot,t)\|_{L^2(\O)}^2\\&\le \|(v-k)_+(\cdot,\si)\|_{L^2(\O)}^2\\&\le C\| (v-k)_+(\cdot,\si)\|_{L^{24/\ve}(\O)}^2|A_k(\si)|^{1-\ve/12}\\&\le C\left(\| (v-k)_+(\cdot,\si)\|_{L^2(\O)}^2+\|\na  (v-k)_+(\cdot,\si)\|_{L^{24/(12+\ve)}(\O)}^2\right)\nu_k^{1-\ve/12}\\&\le C\nu_k^{1+\ve/12},\ea\enn where in the last inequality one has used \eqref{bb12}. This implies
\bnn \la{bb25}\ba  \nu_h  \le C(h-k)^{-2}\nu_k^{1+\ve/12},\ea\enn which, together with the  De Giorgi-type lemma \cite[Lemma 4.1.1]{wyw},
    thus shows \be \la{bb80}\sup_{0\le t\le T}\|\n\|_{L^\infty(\O)}\le \ti C.\ee

Finally,   since $w\triangleq v^{-1}$ satisfies
\be \la{bb81} 2w_t+ 2u\cdot \na w- w\div u+4\ve h_\ve'(\n) w^{-1}|\na w|^2=2\ve\div(h_\ve'(\n)\na w),\ee multiplying \eqref{bb81} by $(w-k)_+$ with $ k\ge \|w(\cdot,0)\|_{L^\infty(\O)} =\|\n_0^{-1/2}\|_{L^\infty(\O)} $ yields that
\be\la{bb82}\ba &\frac{d}{dt}\int (w-k)_+^2dx+2\ve\al \int \n^{\al-1}|\na(w-k)_+ |^2dx \\&\le C\int_{\ti A_k(t)}w|u||\na w|dx+C\int_{\ti A_k(t)}(w-k)_+|u||\na w|dx\\&\le C\int_{\ti A_k(t)}\n^{-\al}|u|^2 dx+\ve \al \int \n^{\al-1}|\na(w-k)_+ |^2dx, \ea\ee where $\ti A_k(t)\triangleq \{x\in\Omega|w(x,t)>k\}.$
It follows from H\"older's inequality, \eqref{e1}, and \eqref{bn1}  that
\be\la{bb83}\ba &\int_{\ti A_k(t)}\n^{-\al}|u|^2 dx \\&=\int_{\ti A_k(t)}\n^{-\al-2/(2+\ve)}(\n^{1/(2+\ve)}|u|)^2 dx\\ &\le C\left(\int_{\ti A_k(t)}\n |u|^{2+\ve}dx\right)^{2/(2+\ve)} \left(\int_{\ti A_k(t)}\n^{-(2+(2+\ve)\al)\ve^{-1}}dx\right)^{\ve/(2+\ve)}\\ &\le C \left(\int_{\ti A_k(t)}\n^{- 6(2+(2+\ve)\al)/(\ve(4-\ve))}dx\right)^{(4-\ve)\ve/(12+6\ve)}|\ti A_k(t)|^{\ve/6}\\ &\le C \ti\nu_k^{\ve/6} ,\ea\ee where $\ti\nu_k\triangleq \sup\limits_{0\le t\le T}|\ti A_k(t)|.$ Hence, putting \eqref{bb83} into \eqref{bb82} leads to
\be\la{bb84}\ba &\frac{d}{dt}\int (w-k)_+^2dx+ \ve\al \int \n^{ \al-1}|\na(w-k)_+ |^2dx \le C \ti\nu_k^{\ve/6}.\ea\ee

Using \eqref{bb84} and \eqref{bn1}, one can proceed in the same way as the proof of \eqref{bb80} to obtain that there exists some positive constant $C\ge \ti C$ such that   \bnn \sup_{(x,t)\in \O\times (0,T)}\n^{-1}(x,t)\le   C,\enn which, combined with \eqref{bb80}, gives \eqref{c1} and finishes the proof of Lemma \ref{lem13}.

We still need the following lemma concerning   the higher order   estimates on $(\n,u)$ which are necessary  to obtain the   global strong solution to the problem  \eqref{ba1}--\eqref{ma1} \eqref{bb35}.
\begin{lemma}\la{lem13'} For any $p> 2,$ there exists some constant $C $ depending on $\ve,p,$ and $T$ such that \be\la{bc4}   \int_0^T\left(\|(\n_t,\na \n_t,u_t\|^p_{L^p(\O) }+\|(\n,\na \n, u)\|_{W^{2,p} (\O) }^p\right)dt\le C . \ee\end{lemma}

  {\it Proof.} First, it follows from \eqref{c1}, \eqref{e1}, \eqref{bb30}, and \eqref{bb5'} that \be \la{bb89}\sup_{0\le t\le T}\left(\|u\|_{L^{2+\ve}(\O)}+\|\na \n\|_{L^2(\O)}\right)+\int_0^T\left(\|u\|^{4+2\ve}_{L^{4+2\ve}(\O)}+\|\na u\|_{L^2(\O)}^2\right)dt\le C,\ee which, together with the standard H\"older estimates for \eqref{bb22},  yields that there exist positive constants $C$ and   $\si\in(0,1)$ such that
\be \la{bb86}\|v\|_{C^{\si,\si/2}(\ol{\O}\times[0,T])}\le C.\ee

Next, it follows from \eqref{bb22} that  $v=\n^{1/2}$ satisfies
\be\la{bb87}2 v_t- 2\ve \div(h_\ve'(\n)\na v ) =-\div(u  v+\na w)- |\O|^{-1}\int u\cdot \na v dx,  \ee where for $t>0,$ $w(\cdot,t)$ is the unique solution to the following problem
\be\la{bb88}\begin{cases} \Delta w=u\cdot \na v- |\O|^{-1}\int u\cdot \na v dx,& x \in \O, \\ \int wdx=0.\end{cases} \ee Applying standard $L^p$-estimates to \eqref{bb88} yields that    $\na w$  satisfies for any $p> 2$
\be \la{bb90} \|\na w\|_{L^p(\O)}\le C(p)\|u\|_{L^{ 2+\ve}(\O)}\|\na v\|_{L^q(\O)}\le C(p)\|\na v\|_{L^q(\O)},\ee where \be\la{c2}\frac{1}{q}\triangleq\frac{1}{p}+\frac{1}{2}-\frac{1}{2+\ve}>\frac{1}{p}. \ee

Since \eqref{bb89} and \eqref{c1}   imply  \bnn \left| \int u\cdot \na v dx \right| \le C\|u\|_{L^2(\O)}\|\na \n\|_{L^2(\O)}\le C,\enn applying  standard parabolic $L^p$-estimates to \eqref{bb87}  yields that for any $p>2 $
\bnn \ba \int_0^T\|\na v \|_{L^p(\O)}^pdt&\le C(p)+C(p)\int_0^T\left(\|u\|_{L^p(\O)}^p+\|\na w\|_{L^p(\O)}^p\right)dt \\ &\le C(p)+C(p)\int_0^T \|u\|_{L^p(\O)}^pdt+\frac{1}{2}\int_0^T \|\na v\|_{L^p(\O)}^p dt,\ea\enn where in the second inequality,   \eqref{bb90},  \eqref{c2}, and \eqref{bb89}  have been used.  Thus,
\be \la{bb91} \int_0^T\|\na v \|_{L^p(\O)}^pdt\le C(p)+C(p)\int_0^T\|u\|_{L^p(\O)}^pdt ,\ee which, together with \eqref{bb89}, gives  \be \la{bb96} \int_0^T\|\na v \|_{L^{4+2\ve}(\O)}^{4+2\ve}dt\le C .\ee

Next, note that \eqref{ba1}$_2$ implies that $u$ satisfies
\be  \la{bb93}   u_t   -   (\frac{ 1}{2}+\sqrt{\ve}) \n^{-1}h_\ve(\n) \Delta u -\left(\frac{ 1}{2} \n^{-1}h_\ve(\n) +   (1+\sqrt{\ve})  \n^{-1}g_\ve(\n)\right) \nabla \div u=F,\ee where \be \la{bg3}\ba F&\triangleq -   u\cdot \na u + (\frac{ 1}{2}+\sqrt{\ve})\n^{-1}\na h_\ve(\n) \cdot \na u+\frac12 \n^{-1}\na u\cdot\na h_\ve(\nv)\\&\quad+    (1+\sqrt{\ve}) \n^{-1} \nabla(g_\ve(\n) )\div u  -\n^{-1}\na P   -e^{-\ve^{-3}}  ( \n^{-1+\ve^{-2}}+\n^{-1-\ve^{-2}}) u.\ea\ee Since $\sqrt{\n}(=v)$ satisfies \eqref{bb86},  applying the standard $L^p$-estimates  to \eqref{bb93} \eqref{bg3}  \eqref{bb35}   with periodic data, we obtain  after using   \eqref{bb96} and \eqref{bb89}     that
  \be  \la{bb99} \ba  \|u_t\|_{L^{2+\ve}(\O\times(0,T))}+ \|\na^2 u \|_{L^{2+\ve}(\O\times(0,T))} \le C+C \|\na u \|_{L^{4+2\ve}(\O\times(0,T))} .  \ea\ee

It thus follows from the Sobolev inequality (\cite[Chapter II (3.15)]{la1}) that for any $\eta>0$ there exists some constant $C(\eta)$ such that \bnn \ba \|\na u \|_{L^{4+2\ve}(\O\times(0,T))}\le & \eta(\|u_t\|_{L^{2+\ve}(\O\times(0,T))}+ \|\na^2 u \|_{L^{2+\ve}(\O\times(0,T))}) \\&+C(\eta)\|u\|_{L^{2+\ve} (\O\times(0,T))} , \ea\enn
which, together with \eqref{bb89} and \eqref{bb99},  gives
 \bnn  \la{bb94} \ba  \|u_t\|_{L^{2+\ve}(\O\times(0,T))}+ \|\na^2 u \|_{L^{2+\ve}(\O\times(0,T))} \le C . \ea\enn 
This, combined with the Sobolev inequality (\cite[Chapter II (3.15)]{la1}), leads to \be  \la{bb95}\sup_{0\le t\le T} \|u\|_{L^\infty(\O)}\le C,\ee which, along with \eqref{bb91}, shows that for any $p>2$
\be \la{bb97} \int_0^T\|\na v \|_{L^p(\O)}^pdt\le C(p).\ee
This, together with \eqref{bb95}, \eqref{c1}, and the   standard $L^p$-estimates of the parabolic system \eqref{bb93}  \eqref{bg3} \eqref{bb35},   yields that for any $p>4, $
 \bnn\ba &\|u_t\|_{L^p(\O\times(0,T))}+ \|\na^2 u \|_{L^p(\O\times(0,T))}\\ &\le C(p)+C (p)\|\na u \|_{L^2(\O\times(0,T))}^{1/p}\|\na u \|_{L^\infty(\O\times(0,T))}^{1-1/p} \|\na \n \|_{L^{2p}(\O\times(0,T))}\\&\le C(p)+\frac{1}{2}\|u_t\|_{L^p(\O\times(0,T))}+\frac{1}{2} \|\na^2 u \|_{L^p(\O\times(0,T))} ,\ea\enn where in the second inequality one has used \eqref{bb89} and  the Sobolev inequality (\cite[Chapter II (3.15)]{la1}).   Thus, it holds that for any $p>2 , $
 \be \la{c3} \|u_t\|_{L^p(\O\times(0,T))}+ \|\na^2 u \|_{L^p(\O\times(0,T))} \le C(p).\ee
With \eqref{c3} and  \eqref{bb97} at hand, one can deduce easily from \eqref{bb22}, \eqref{bb35}, and \eqref{bb86} that for any $p> 2,$
 \bnn \|\n_t \|_{L^p( 0,T,W^{1,p}(\O))}+ \|\na^2 \n \|_{L^p( 0,T,W^{1,p}(\O))} \le C(p),\enn which, together with \eqref{c3} and \eqref{bb89},  gives the desired estimate \eqref{bc4} and finishes the proof of Lemma \ref{lem13'}.


\subsection{\la{qsec1}Compactness results}

 Throughout this subsection, it will be  always assumed that $\al$ and $\ga$ satisfy the conditions listed in Theorem \ref{th2}.

We first construct the  initial data. Set \be\la{lk1}\si_0\triangleq(8(\al+\ga+2))^{-8}.\ee
Choose
\be\nonumber\la{pd5}\ti \n_{0\ve}\in C^\infty (\O),\quad  0 \le\ti\n_{0\ve} \le \ve^{-4\si_0}   \ee  satisfying \be\nonumber\la{pd17} \| \ti\n_{0\ve}-\n_0\|_{L^1 (\O)}+ \| \ti\n_{0\ve}-\n_0\|_{L^\ga (\O)}+\|\na ( \ti\n_{0\ve}^{\al-1/2}-\n_{0}^{\al-1/2})\|_{L^2(\O)}<\ve.\ee For $\nu\ge 2$  suitably large such that
$ \nu(\al-1/2)\ge 5,$ define \be\la{pd7}\n_{0\ve}=\left(\ti\n_{0\ve}^{\nu(\al-1/2)}+\ve^{4\si_0 \nu (\al-1/2)}\right)^{2/(\nu(2\al-1))}.\ee It is easy to check that
\be\la{pd8} \lim_{\ve\rightarrow 0}\|\n_{0\ve}  -\n_0\|_{L^1(\O)}=0\ee and that there exists some constant $C$ independent of $\ve$ such that \eqref{pd9} holds.

Since $\n_0,m_0$ satisfy \eqref{pini1},   we choose $w_{0\ve} \in C^\infty(\O)$ such that \be\nonumber \|w_{0\ve} - m_0/\n_0^{(1+\eta_0)/ (2+\eta_0)}\|_{L^{2+\eta_0}(\O)}\le \ve .\ee
Set   \be \la{pd21} u_{0\ve}=  \n_{0\ve}^{-1/(2+\eta_0)}w_{0\ve} .\ee  Then, we have
\be\la{pd13}\lim_{\ve\rightarrow 0}\|\n_{0\ve} u_{0\ve} -m_0\|_{L^1(\O)}=0.\ee  Moreover, there exists some positive constant $C$ independent of $\ve$ such that \eqref{pd22} holds.

Extend then  $(\n_{0\ve},u_{0\ve} )$ $\O$-periodically to $\r^2.$   The standard   parabolic  theory  \cite{la1}, together with Lemmas \ref{lem13} and \ref{lem13'}, thus yields that
the problem \eqref{ba1}--\eqref{ma1} \eqref{bb35},   where the initial data  $(\n_{0 },u_{0 })$ is replaced by $(\n_{0\ve},u_{0\ve}),$   has  a unique strong solution   $( \n_\ve,u_\ve) $    satisfying \bnn \nv,\,\,\uv , \,\, \n_{\ve t} ,\,\,u_{\ve t}  ,\,\, \na^2\nv,\,\,\na^2\uv \in L^p(\O\times (0,T)),\enn for any $T>0$ and  any $p>2.$ Moreover,  all estimates obtained by Lemmas \ref{lem10}-\ref{lem01} still hold for $(\nv,\uv).$

Letting $\ve\to 0^+,$ we will modify  the compactness results due to \cite{MelletVasseur05} to prove that   the limit (in some sense) $(\n ,\sqrt{\n }u)$ of $(\nv,\sqrt{\nv}\uv) $ (up to a subsequence) is a weak solution to
\eqref{ii1} \eqref{aa'1} \eqref{hgvv1} \eqref{en1}. We begin with the following strong convergence of $\nv.$
\begin{lemma} \la{lema2}   There exists a function $\n\in L^\infty(0,T;L^1(\O)\cap L^\ga(\O))$ such that up to a subsequence,
 \be \la{e3}\nv \ro \n  \mbox{   in }L^\ga (\O\times(0,T )).\ee In particular, \be \la{e3'}\nv  \ro \n  \mbox{ almost everywhere  in  } \O\times(0,T  ).\ee \end{lemma}
  {\it Proof.} First, it follows from \eqref{bb30} and  \eqref{bb5'}  that there exists some generic positive constant $C$ independent of $\ve  $  and $T  $     such that
\be\la{d1}\ba& \sup_{0\le t\le T}\int  (\n_\ve |u_\ve|^2+\n_\ve+\n_\ve^\ga)dx    +\ioo h_\ve(\n_\ve)|\na u_\ve|^2 dxdt \\&+\ve\ioo \nv^{-1} h_\ve'(\nv) |\na \nv|^2(1+|\uv|^2)   dx  dt \\&+e^{-\ve^{-3}}\ioo \left(      \n_\ve^{ \ve^{-2}}+ \n_\ve^{-\ve^{-2}}\right) |u_\ve|^2dx  dt  \le C ,\ea\ee
and that
\be\la{d2}\ba& \sup_{0\le t\le T} \int  \nv^{-1}(h_\ve'(\nv))^2|\na\nv|^2 dx   + \ioo  \n_\ve^{\al+\ga-3}   |\na \n_\ve|^2dxdt \\&+\ve^{13/3}e^{-\ve^{-3}}\sup_{0\le t\le T}\int   \left( \nv^{ \ve^{-2}+ \ti\ga-1 }+ \nv^{-\ve^{-2}-1/8} \right)  dx  \le C.\ea\ee

Then,  \eqref{d1} and \eqref{d2} imply that
\be \la{e4}\sup_{0\le t\le T}\|\na \nv^\al\|_{L^{2\ga/(\ga+1)}(\O)  }\le C\sup_{0\le t\le T}\|\nv^{1/2}\|_{L^{2\ga} (\O)  }\sup_{0\le t\le T} \|\na \nv^{\al-1/2}\|_{L^2 (\O)  }\le C.\ee Moreover, note that $\nv^\al$ satisfies
\be \la{e5}\ba& (\n_\ve^\al)_t+\div (\n_\ve^\al u_\ve)+(\al-1)\n_\ve^{\al}\div u_\ve \\&=\ve\al\div(\nv^{\al-1}h_\ve'(\nv)\na\nv)-\ve\al(\al-\frac{1}{2})\nv^{\al-2}h_\ve'(\nv)|\na \nv|^2 . \ea\ee It follows from \eqref{d1},   \eqref{d2}, and \eqref{bb6'} that
\be \la{e6}\ba  \sup_{0\le t\le T}\|\n_\ve^\al u_\ve\|_{L^1 (\O)  }  \le C\sup_{0\le t\le T}\|\nv^{\al-1/2}\|_{L^2 (\O)  }\sup_{0\le t\le T}\| \nv^{1/2}\uv\|_{L^2 (\O)  } \le C ,\ea\ee
  \be \la{e7'}\ba  \int_0^T\|\n_\ve^{\al}\div u_\ve\|^2_{L^{1} (\O)  }dt  \le C\int_0^T\|\nv^{\al/2}\|^2_{L^2  (\O) }\|\nv^{\al/2}\div u_\ve\|^2_{L^{2} (\O)  }dt
\le C ,\ea\ee and that \be \la{e7}\ba & \sup_{0\le t\le T}\int \left(\nv^{\al-1}h_\ve'(\nv)|\na\nv|+\nv^{\al-2}h_\ve'(\nv)|\na \nv|^2\right)dx
\\&\le C \sup_{0\le t\le T}\int\left(\nv^{2\al-1} +\nv^{-1}(h_\ve'(\nv))^2|\na \nv|^2\right)dx  \le C .\ea\ee
The combination of   \eqref{e5}--\eqref{e7} implies that \be\la{ve4}\|(\nv^\al)_t\|_{L^2(0,T;W^{-1,1} (\O) )}\le C.\ee Letting $\ve\to 0^+,$ it follows from \eqref{e4}, \eqref{ve4}, and the
Aubin-Lions lemma  that  up to a subsequence \bnn\la{e2} \nv^\al \ro \n^\al \mbox{ in }C([0,T];L^{3/2} (\O)),\enn   which
 implies that \eqref{e3'} holds. In   particular, it holds that \be\la{e8} \nv^{\al-1/2} \ro \n^{\al-1/2} \mbox{ in }L^2(0,T ;L^2   (\O)).\ee

Finally, it follows from the Sobolev inequality, \eqref{d1}, and \eqref{d2} that
\be\la{e9}\ba &\int_0^T\|\nv^\ga\|^{(5\ga+3(\al-1))/(3\ga)} _{L^{(5\ga+3(\al-1))/(3\ga)}(\O)   }dt
\\& \le C\int_0^T\|\nv\|_{L^\ga (\O)  }^{2\ga/3}\left(\|\nv\|_{L^1(\O) }^{\ga+\al-1}+\|\na \nv^{(\ga+\al-1)/2}\|_{L^2 (\O)  }^2\right)dt\le C,\ea\ee which together with \eqref{e3'} thus gives \eqref{e3} due to  $(5\ga+3(\al-1))/(3\ga)>1.$ The proof of Lemma \ref{lema2} is finished.

Before proving the strong convergence of $ \sqrt{\nv}\uv $ in $L^2(\O\times (0,T)),$ we show first the following compactness   of $\nv^{(\ga+1)/2}\uv.$

\begin{lemma} \la{lem02}  There exists a function $m(x,t)\in L^2(\O\times(0,T ))$ such that up to a subsequence,
\be \la{e25}\nv^{(\ga+1)/2}\uv\rt m \mbox{ in }L^2(0,T;L^p (\O)) ,\ee      for all $p\in [1,2).$
Moreover, \be \la{e25'}\nv^{(\ga+1)/2}\uv\rt \n^{(\ga+1)/2}u\mbox{ almost everywhere }(x,t)\in\O\times (0,T) ,\ee  where \be \la{e24}u(x,t)\triangleq\begin{cases} m(x,t)/\n^{(\ga+1)/2}(x,t)& \mbox{ for } \n(x,t)>0,\\0,&\mbox{ for } \n(x,t)=0.\end{cases}\ee
 \end{lemma}
  {\it Proof.}   First, since  $\al\in (1/2,(1+\ga)/2],$ it follows from
     \eqref{d1},   \eqref{d2}, and   \eqref{bb6'}  that for any $\eta>0,$
\be\la{e27}\ba &\int_0^T\|\na (\nv^{(1+\ga)/2} \uv)\|_{L^1  (\O) }^2dt\\&\le  C\int_0^T\|\nv^{(1+\ga-\al)/2}\|_{L^2   (\O)}^2 \|\nv^{\al/2}\na \uv\|_{L^2  (\O) }^2dt\\&\quad+ C\int_0^T\|\nv^{(\ga-2\al+2)/2}\uv\|_{L^2(\O)   }^2   \|\na \nv^{ \al-1 /2}\|_{L^2(\O)  }^2   dt \\&\le C(\eta)+C(\eta)\int_0^T\|  \nv^{1/2} \uv \|_{L^2 (\O)   }^2dt+C\eta\int_0^T\|  \nv^{(1+\ga)/2} \uv \|_{L^2(\O)   }^2dt   , \ea\ee which together with the Sobolev  inequality gives
\be\la{e29}\int_0^T\|\na (\nv^{(1+\ga)/2} \uv)\|_{L^1 (\O)  }^2dt+\int_0^T\|  \nv^{(1+\ga)/2} \uv \|_{L^2  (\O) }^2dt\le C,\ee
due to the following simple fact
\be\la{e30} \sup_{0\le t\le T}\int  \nv^{(\ga+1)/2}|\uv|dx\le C\sup_{0\le t\le T}\|\nv^\ga\|_{L^1  (\O)  }^{1/2}\sup_{0\le t\le T}\|\nv^{1/2}\uv\|_{L^2 (\O)  }\le C.\ee

Next,  we claim that \be \la{ux1}\|(\nv^{(\ga+1)/2}\uv)_t\|_{L^1(0,T;W^{-1,1} (\O))}\le C,\ee
which, combined with \eqref{e29} and the Aubin-Lions lemma, yields that  there exists a function $m(x,t)\in L^2(\O\times(0,T ))$ such that up to a subsequence, \eqref{e25} holds   for all $p\in [1,2).$
In particular, \be \la{e22}\nv^{(\ga+1)/2}\uv\rt m \mbox{ almost everywhere }(x,t)\in\O\times (0,T)  .\ee
Moreover, since $\nv^{1/2}\uv$ is bounded in $L^\infty(0,T;L^2(\O)),$ Fatou's lemma gives
\be\nonumber \la{e23}\ioo  \liminf\frac{|\nv^{(\ga+1)/2}\uv|^2}{\nv^\ga}dxdt<\infty,\ee which implies  $m(x,t)=0$ almost everywhere in $\{(x,t)\in\O\times(0,T)|\n(x,t)=0\}.$  Hence, for $u(x,t)$ as in \eqref{e24},   we arrive at \be\nonumber\la{e26} m(x,t)=\n^{(\ga+1)/2}(x,t)u(x,t),\ee    which together with  \eqref{e22} gives   \eqref{e25'}.

Finally, it remains to prove  \eqref{ux1}. In fact,  note that
\be\la{e31} \ba (\nv^{(\ga+1)/2} \uv)_t&=\frac{\ga+1}{2}\nv^{(\ga-1)/2}(\nv)_t \uv+\nv^{(\ga+1)/2}   (\uv)_t  . \ea\ee
One can use \eqref{ba1} to get
\be\la{e32} \ba  \nv^{(\ga-1)/2}(\nv)_t \uv &=-\nv^{(\ga-1)/2}\div (\nv\uv)\uv+\ve\nv^{\ga/2}\div (\nv^{-1/2}h'_\ve(\nv)\na\nv)\uv\\&=-\frac{2}{\ga+1}\div(\nv^{(\ga+1)/2}\uv\otimes \uv)-\frac{\ga-1}{\ga+1}\nv^{(\ga+1)/2}\div\uv \uv\\&\quad+\frac{2}{\ga+1}\nv^{(\ga+1)/2} \uv \cdot\na\uv+ \ve\div (\nv^{(\ga-1)/2}h'_\ve(\nv)\na\nv\otimes\uv)\\&\quad-\frac{\ga\ve}{2}\nv^{(\ga-3)/2} h'_\ve(\nv)|\na \nv|^2\uv-\ve \nv^{(\ga-1)/2}h'_\ve(\nv)\na\nv\cdot\na\uv,\ea\ee
and
\be\la{e33} \ba  \nv^{(\ga+1)/2}   (\uv)_t =&- \nv^{(\ga+1)/2}\uv\cdot\na \uv +\div(\nv^{(\ga-1)/2}h_\ve(\nv)(\mathcal{D}u+\sqrt{\ve}\na\uv))\\&-\frac{\ga-1}{2}\nv^{(\ga-3)/2} h_\ve(\nv)\na\nv\cdot(\mathcal{D}u+\sqrt{\ve}\na\uv)) \\&+(1+\sqrt{\ve})\na(\nv^{(\ga-1)/2}g_\ve(\nv)\div\uv)\\&-\frac{(\ga-1)(1+\sqrt{\ve})}{2}\nv^{(\ga-3)/2} g_\ve(\nv)\na\nv\div\uv  \\&-\nv^{(\ga-1)/2}\na\nv^\ga-e^{-\ve^{-3}} ( \nv^{ \ve^{-2}+(\ga-1)/2} +\nv^{-\ve^{-2}+(\ga-1)/2} )\uv. \ea\ee

One needs to estimate each term on the righthand side of  \eqref{e32} and \eqref{e33}. It follows from the H\"{o}lder inequality, \eqref{d1},   \eqref{d2},  \eqref{e29}, and the Sobolev inequality that
\be\la{e34} \ba \int_0^T\|\nv^{(\ga+1)/2}|\uv|^2\|_{L^1  (\O) }dt \le  C\int_0^T\|(\nv +\nv^{ \ga+1 })|\uv|^2\|_{L^1(\O)   } dt \le C,\ea\ee
\be\la{e35} \ba& \int_0^T\|\nv^{(\ga+1)/2}|\uv||\na\uv|\|_{L^1 (\O)  }dt \\&\le  C\int_0^T\|(\nv^{1/2} +\nv^{ (\ga+1)/2 })|\uv| \|_{L^2 (\O)   }^2dt+C\int_0^T\|\nv^{\al/2}|\na \uv|\|_{L^2  (\O)}^2 dt \le C,\ea\ee
\be\la{e36} \ba& \ve\int_0^T\|\nv^{(\ga-1)/2}h'_\ve(\nv)|\na\nv||\uv|\|_{L^1(\O)   }dt \\&\le  C\ve\ioo  \nv^{-1} h_\ve'(\nv)|\na\nv|^2 |\uv|^2 dxdt  + C \ioo \nv^{\ga-1} h_\ve(\nv) dx dt\\&\le C+C\ioo (\nv^{\ga+\al-1}+\nv^{\ga -1/8}+\nv^{ 2\ga-5/6})dt \le C ,\ea\ee
\be\la{e37} \ba& \ve\int_0^T\|\nv^{(\ga-3)/2} h'_\ve(\nv)|\na \nv|^2|\uv|\|_{L^1 (\O) }dt \\&\le  C\ve\ioo \left( \nv^{-1} h_\ve'(\nv)|\na\nv|^2|\uv|^2    +  \n^{\ga-3}h_\ve(\nv)|\nabla \nv|^2\right)dxdt \le C,\ea\ee
 \be\la{e38} \ba&  \ioo \nv^{(\ga-1)/2}(h_\ve(\nv)+|g_\ve(\nv)|)| \na\uv |dxdt   \\&\le  C \ioo h_\ve(\nv)|\na\uv|^2  dxdt +C \ioo  \nv^{\ga-1}h_\ve(\nv) dxdt \le C,\ea\ee\be\la{e39} \ba&  \ioo \left(  \nv^{(\ga-1)/2}h'_\ve(\nv)+  \nv^{(\ga-3)/2}(h_\ve(\nv)+|g_\ve(\nv)|)\right) |\na\nv||\na\uv|dx dt \\&\le C\ioo  \nv^{(\ga-3)/2} h_\ve(\nv) |\na\nv||\na\uv|dx dt \\&\le  C \ioo   h_\ve(\nv)|\na\uv|^2  dxdt +C \ioo \nv^{\ga-3}h_\ve(\nv)|\nabla \nv|^2dxdt \le C,\ea\ee
 \be\ba &\ioo \nv^{(\ga-1)/2}|\na\nv^\ga|dxdt\\&\le C\int_0^T\|\nv^{\ga-\al/2}\|_{L^2  (\O) }\|\na\nv^{(\ga+\al-1)/2}\|_{L^2  (\O) }dt\le C ,\ea\ee

\be\la{e40} \ba&  e^{-\ve^{-3}}\ioo  (\nv^{ \ve^{-2}+(\ga-1)/2} + \nv^{-\ve^{-2}+(\ga-1)/2})|\uv|dxdt\\& \le C e^{-\ve^{-3}}\ioo ( \nv^{ \ve^{-2} }+\nv^{-\ve^{-2} })|\uv|^2dxdt\\&\quad+ C e^{-\ve^{-3}}\ioo  (\nv^{ \ve^{-2}+\ga-1 }+\nv^{-\ve^{-2}+\ga-1 }) dxdt\\& \le C + C e^{-\ve^{-3}}\ioo  (\nv^{ \ve^{-2}+\ga-1 }+\nv^{-\ve^{-2} }) dxdt \le C , \ea\ee
 where in the last inequality one has used the following simple facts that
\be \la{e19'}\ba& e^{-\ve^{-3}}\ioo  \nv^{ \ve^{-2}+\ga-1}  dxdt\\&= \frac{e^{-\ve^{-1} /  (6(1+(\ti\ga-1)\ve^2) )}}{\ve^{13(1+(\ga-1)\ve^2)/(3(1+(\ti\ga-1) \ve^2) )}}\ioo  (\ve^{13/3}e^{-\ve^{-3}}\nv^{\ve^{-2}+\ti\ga-1})^{\frac{\ve^{-2}+\ga-1} {\ve^{-2}+\ti\ga-1}} dxdt\\&\le \frac{C  e^{-1 /(9\ve)}}{\ve^{13}}\sup_{0\le t\le T}\left(\ve^{13/3}e^{-\ve^{-3}}\int   \nv^{ \ve^{-2}+\ti\ga -1}   dx\right)^{\frac{\ve^{-2}+\ga-1} {\ve^{-2}+\ti\ga-1}}\\&\le    {C  \ve^{- 13}}{e^{-1 /(9\ve)}}\rightarrow 0\quad \mbox{ as } \ve\rightarrow 0,\ea\ee and that
\be \la{e19}\ba& e^{-\ve^{-3}}\ioo  \n^{-\ve^{-2}}  dxdt\\&= \frac{e^{-\ve^{-1}/( 8+\ve^2) }}{\ve^{104/(3(8+\ve^2))}}\ioo  \left(\ve^{13/3}e^{-\ve^{-3}}   \n^{-\ve^{-2}-1/8}   \right)^{8/(8+\ve^2)} dxdt\\&\le \frac{C  e^{-1/(9\ve )}}{\ve^{5}}\sup_{0\le t\le T}\left(\ve^{13/3}e^{-\ve^{-3}}\int   \n^{-\ve^{-2}-1/8}   dx\right)^{8/(8+\ve^2)}\\&\le    {C  \ve^{- 5}}{e^{-1/(9\ve )}}\rightarrow 0\quad \mbox{ as } \ve\rightarrow 0.\ea\ee
Thus, all these estimates  \eqref{e34}--\eqref{e40} together with \eqref{e31}-\eqref{e33}  yield \eqref{ux1}.
 The proof of Lemma \ref{lem02} is completed.

Now we are in a position to prove the strong convergence of $ \sqrt{\nv}\uv $ in $L^2(\O\times (0,T))$ which in fact is essential to obtain the existence of global weak solution to the problem \eqref{ii1} \eqref{aa'1} \eqref{hgvv1} \eqref{en1}.
\begin{lemma}\la{lema1}   Up to a subsequence, \be \la{e42} \sqrt{\nv}\uv\rt \sqrt{\n}u   \mbox{ strongly in } L^2 ( 0,T;L^2  (\O)),\ee with \be \la{e'42} \sqrt{\n}u \in L^\infty(0,T;L^2(\O )).\ee\end{lemma}
{\it Proof. }
  First, Lemma \ref{lem01} yields that there exists
some constant $C$ independent of $\ve$ such that
\be\la{d3}\sup_{0\le t\le T}\int \n_\ve  |u_\ve|^2 \ln (e+|u_\ve|^2)dx\le C,\ee  which, together with  \eqref{e25'},   \eqref{e3'}, and Fatou's lemma,  gives  \be\la{e44} \ba \ioo \n|u|^2\ln (e+|u|^2)dxdt&\le \ioo \liminf_{\ve\rt 0} \nv|\uv|^2\ln(e+|\uv|^2)dxdt\\&\le \liminf_{\ve\rt 0} \ioo  \nv|\uv|^2\ln(e+|\uv|^2)dxdt\le C.\ea\ee

Next,
 direct calculation shows that for any $M>0,$\be \la{e47} \ba &\ioo   |\sqrt{\nv}\uv-\sqrt{\n}u|^2dxdt\\&\le 2\ioo   |\sqrt{\nv}\uv1_{(|\uv|\le M)}-\sqrt{\n}u1_{(|u|\le M)}|^2dxdt\\&\quad+2\ioo  |\sqrt{\nv}\uv1_{(|\uv|\ge M)}|^2dxdt+2\ioo  |\sqrt{\n }u1_{(|u|\ge M)}|^2dxdt.\ea\ee

Next, it follows  from \eqref{e25'}  and  \eqref{e3'}  that $ \sqrt{\nv}\uv  $ converges almost everywhere to $\sqrt{\n}u  $ in the set $\{(x,t)\in \O\times(0,T)|\n(x,t)>0\}  . $   Moreover, since \be \la{e3a} \sqrt{\nv}|\uv|1_{(|\uv|\le M)}\le M\sqrt{\nv}  ,\ee and  $\nv\rt 0$ almost everywhere in the set $\{(x,t)\in \O\times(0,T)|\n(x,t)=0\}  , $
we have \be \nonumber\la{e45} \sqrt{\nv}\uv1_{(|\uv|\le M)} \rt \sqrt{\n}u1_{(| u|\le M)} \mbox{ almost everywhere in } \O\times(0,T), \ee which, together with \eqref{e3a} and  \eqref{e3},  implies   \be \la{e48} \ioo  \left|\sqrt{\nv}\uv1_{(|\uv|\le M)}-\sqrt{\n}u1_{(|u|\le M)}\right|^2dxdt \rt 0.\ee

 Next, it follows from \eqref{d3} and  \eqref{e44} that  \be \la{e50} \ba & \ioo  |\sqrt{\nv}\uv1_{(|\uv|\ge M)}|^2dxdt+\ioo  |\sqrt{\n }u1_{(|u|\ge M)}|^2dxdt\\&\le \frac{1}{\ln (e+M^2)}\ioo   \left( \nv |\uv|^2\ln(e+|\uv|^2)+ \n |u|^2\ln(e+|u|^2)\right)dxdt\\&\le \frac{C}{\ln (e+M^2)}.\ea\ee
Substituting  \eqref{e48} and \eqref{e50} into  \eqref{e47} yields that up to a subsequence \be \la{e51} \limsup_{\ve\rt 0}\ioo   |\sqrt{\nv}\uv-\sqrt{\n}u|^2dxdt\le \frac{C}{\ln(e+M^2)}\ee
 for any $M>0.$ We thus obtain  \eqref{e42}  by taking $M\rt \infty$ in \eqref{e51}.

  Finally, the combination of \eqref{d1} with  \eqref{e42}  gives   \eqref{e'42} immediately.  The proof of Lemma \ref{lema1} is completed.

 As a consequence of Lemmas \ref{lema2} and  \ref{lema1}, the following convergence of the diffusion terms holds.
\begin{lemma}\la{plem1} Up to a subsequence, \be \la{e54} \nv^\al \na\uv\rt \n^\al\na u \mbox{ in }\mathcal{D}',\ee\be \la{e55} \nv^\al (\na\uv)^{\rm tr}\rt \n^\al(\na u)^{\rm tr}\mbox{ in }\mathcal{D}',\ee\be \la{e56} \nv^\al \div\uv\rt \n^\al\div u \mbox{ in }\mathcal{D}'.\ee\end{lemma} {\it Proof. }
 Let $\phi$ be a test function. Then it follows from \eqref{d2}, \eqref{e8}, and \eqref{e42} that  \bnn \la{e57} \ba &\ioo  \nv^\al \na\uv \phi dxdt\\&=-\ioo  \nv^{\al-1/2}\sqrt{ \nv}  \uv \na\phi dxdt-\frac{2\al}{2\al-1}\ioo  \na \nv^{\al-1/2} \sqrt{\nv}  \uv \phi dxdt\\&\rt-\ioo  \n^{\al-1/2} \sqrt{\n}   u\na\phi dxdt-\frac{2\al}{2\al-1}\ioo  \na \n^{\al-1/2} \sqrt{\n}   u \phi dxdt,\ea\enn
which gives \eqref{e54}.
Similar arguments prove  \eqref{e55} and \eqref{e56}, and   the proof of Lemma \ref{plem1} is completed.

\subsection{Proof  of Theorem  \ref{th2}: $\O=\mathbb{T}^2$}
First,  rewrite  \eqref{ba1}$_1$ as
\be \la{e5'}\ba& (\nv )_t+\div (\n_\ve  u_\ve) =\ve \div( h_\ve'(\nv)\na\nv)-\frac{\ve}{2} \nv^{-1}h_\ve'(\nv)|\na \nv|^2 . \ea\ee
It follows from \eqref{d1} and \eqref{d2} that
\be \la{e135'}  \ba& {\ve} \ioo  h_\ve'(\nv) |\na \nv| dxdt\\&\le C \ve\left( \ioo  \nv^{-1}(h_\ve'(\nv))^2|\na \nv|^2 dxdt\right)^{1/2}\\&\le C {\ve} ,\ea \ee and that
\be \la{e153} \ba&  \ve\ioo  \nv^{-1} h_\ve'(\nv)  |\na \nv|^2 dxdt \\&\le C {\ve}  \left(  \ioo  \nv^{-1} (h_\ve'(\nv))^2  |\na \nv|^2 dxdt\right)^{1/2}  \left(  \ioo  \nv^{-1}   |\na \nv|^2 dxdt\right)^{1/2}\\&\le C\ve^{2/3}\left( \ve^{2/3}\ioo  \left( \nv^{-5/4}+\nv^{2\ga-8/3}\right)   |\na \nv|^2 dxdt\right)^{1/2} \\&\le C\ve^{2/3}.\ea  \ee

Then, letting $\psi$ be a test function, multiplying \eqref{e5'} by $\psi$, integrating the resulting equality over $\O\times (0,T),$ and taking $\ve\rt 0$ (up to a subsequence), one can verify easily after using \eqref{e3}, \eqref{e42},  \eqref{pd8}, \eqref{e135'}, and \eqref{e153} that $(\n,\sqrt{\n}u)$ satisfies \eqref{fin1}.

Next,    \eqref{ba1} implies that
\be\la{e'14}\ba &(\nv\uv)_t+\div(\nv\uv\otimes\uv)-\div(\nv^\al\mathcal{D}\uv) -(\al-1)\na(\nv^\al\div\uv)+\na P(\nv)\\&=\ve\div(h'_\ve(\nv)\na\nv\otimes \uv)-\frac{\ve}{2}\nv^{-1}h'_\ve(\nv)|\na\nv|^2\uv-\ve h'_\ve(\nv)\na\nv \cdot \na\uv\\&\quad-e^{-\ve^{-3}}(\nv^{ \ve^{-2}}+\nv^{-\ve^{-2}}) \uv+\sqrt{\ve}\div(h_\ve(\nv) \na\uv ) +\sqrt{\ve}\na(g_\ve(\nv) \div\uv) \\&\quad+\ve^{1/3}\div((\nv^{7/8}+\nv^{\ti\ga})\mathcal{D}\uv) -\frac{\ve^{1/3}}{8}\na(\nv^{7/8}\div\uv)+\ve^{1/3}(\ti\ga-1)\na(\nv^{\ti\ga}\div\uv).\ea\ee

  Using \eqref{d1} and \eqref{d2}, we have
 \be \la{e152} \ba& {\ve} \ioo  h_\ve'(\nv) |\na \nv| |\uv|dxdt\\&\le C {\ve} \left( \ioo   \nv^{-1} (h_\ve' (\nv) )^2|\na \nv|^2 dxdt  \right)^{1/2}\left( \ioo \nv |\uv|^2 dxdt\right)^{1/2}\\&\le C   {\ve}\rt 0,\ea \ee
 \be \la{e153'} \ba& {\ve} \ioo \nv^{-1} h_\ve'(\nv) |\na \nv|^2 |\uv|dxdt\\&\le C   \left(\ve\ioo   \nv^{-1} h_\ve' (\nv) |\na \nv|^2 |\uv|^2dxdt  \right)^{1/2}\\&\quad\times \left( \ve\ioo \nv^{-1} h_\ve'(\nv)  |\na \nv|^2 dxdt\right)^{1/2} \\&\le C\ve^{1/3}\rt 0,\ea  \ee where \eqref{e153} has been used  in the last inequality,
 \be \ba& {\ve} \ioo   h_\ve'(\nv) |\na \nv| |\na\uv|dxdt\\&\le C \sqrt{\ve} \left(\ve\ioo   \nv^{-1} h_\ve' (\nv) |\na \nv|^2  dxdt  \right)^{1/2}\left( \ioo h_\ve (\nv)|\na\uv|^2 dxdt\right)^{1/2}\\&\le C \sqrt{\ve}\rt 0,\ea \ee
 \be\la{e52'}\ba &e^{-\ve^{-3}}\ioo (\nv^{ \ve^{-2}}+\nv^{-\ve^{-2}})| \uv|dxdt\\&\le C \left(e^{-\ve^{-3}}\ioo (\nv^{ \ve^{-2}}+\nv^{-\ve^{-2}})| \uv|^2dxdt\right)^{1/2} \\&\quad\times  \left(e^{-\ve^{-3}}\ioo (\nv^{ \ve^{-2}}+\nv^{-\ve^{-2}})  dxdt\right)^{1/2}\\&\le C\left(e^{-\ve^{-3}}\ioo (\nv +\nv^{ \ve^{-2}+\ga-1}+\nv^{-\ve^{-2}})  dxdt\right)^{1/2}\rt 0   \ea\ee due to  \eqref{e19'} and  \eqref{e19},
\be \la{e52} \ba&\sqrt{\ve} \ioo (h_\ve(\nv)+|g_\ve(\nv)|)|\na \uv| dxdt\\&\le C\sqrt{\ve}\left( \ioo  h_\ve(\nv)  dxdt\right)^{1/2}\left( \ioo  h_\ve(\nv)|\na \uv|^2 dxdt\right)^{1/2}\\&\le C \sqrt{\ve}\rt 0,\ea \ee where  in the second inequality   one has used the fact that
\be \la{e5'3}\ba \int_0^T \|\nv\|^{\ti\ga }_{L^{\ti\ga }(\O)  }dt&=\int_0^T\|\nv\|^{ \ga+1/6}_{L^{ \ga+1/6}  (\O)} dt \\&\le  C  \int_0^T\left( \|\nv \|_{L^1 (\O) }+   \| \nv \|^{5\ga/3+\al-1} _{L^{5\ga/3+\al-1} (\O) }\right)dt\\&\le  C  ,\ea\ee due to $\ga>1,$ $\al>1/2,$ and \eqref{e9},
 \be \la{e53} \ba& \ve^{1/3} \ioo    (\nv^{7/8}+\nv^{\ti\ga })|\na \uv| dxdt \\&\le C\ve^{1/6}\left( \ioo  (\nv^{7/8}+\nv^{\ti\ga }) dxdt\right)^{1/2}\left( \ioo  h_\ve(\nv)|\na \uv|^2 dxdt\right)^{1/2}\\&\le C \ve^{1/6} \rt 0, \ea \ee where \eqref{e5'3} has been used in   the second inequality.

Finally, let $\phi$ be a test function. Multiplying \eqref{e'14} by $\phi,$   integrating the resulting equality over $\O\times (0,T),$  and taking $\ve\rt 0$ (up to a subsequence),  by Lemmas  \ref{lema2},  \ref{lema1},  and \ref{plem1},  we obtain after using \eqref{e152}--\eqref{e52}, \eqref{e53},  and \eqref{pd13} that $(\n,\sqrt{\n}u)$ satisfies \eqref{fin2}. The proof  of Theorem  \ref{th2} in the case $\O=\mathbb{T}^2$ is completed.


\section{Proof of Theorem  \ref{vth2}: $\O=\mathbb{T}^3$}

In this section, we will show how to modify the analysis in the previous section to deal with the 3-dimensional case with periodic boundary conditions.

\subsection{A priori estimates}
   For $\ve$ as in \eqref{ma1}, let $\O=\mathbb{T}^3 $
and  smooth  functions $ \n_{0\ve}>0$ and $u_{0\ve} $ satisfy
 \eqref{pd9}  and \eqref{pd22} for some constant $C$ independent of $\ve.$ Moreover, if $\al\in (1,2),$ in addition to \eqref{pd9}  and \eqref{pd22},  it holds that for some $C$ independent of $\ve,$ \be \la{pd3q} \int \n_{0\ve}|u_{0\ve}|^4dx\le C.\ee
We extend  $\n_{0\ve} $ and $u_{0\ve}$ $\O$-periodically to $\r^3. $  For  $\ti\ga$ as in \eqref{ma1},   consider
     \be\la{vba1} \begin{cases}   \n_t+\div(\n u)  =  \ve \n^{1/2}\div(\n^{- 1/2}h_\ve'(\n)\na \n),
     \\   \n u_t+ \n u\cdot \na u   -\div(h_\ve(\n)   \na u)   -\na(g_\ve(\n)\div u)  +\na P\\ =-e^{-\ve^{-3}}  (\n^{ \ve^{-2}}+ \n^{-\ve^{-2}})u ,   \end{cases}\ee
with  $h_\ve(\n)$ and $g_\ve(\n)$  as in \eqref{ini9}.    The initial   condition  for the system \eqref{vba1} is   imposed as:
\be \la{vbb35'}  (\n,u)(x,0)=(\n_{0\ve},u_{0\ve}),\quad x\in \O. \ee
Let $T>0$ be a fixed time and $(\rho,u)$ be
the smooth solution to \eqref{vba1} \eqref{ini9} \eqref{ma1}  \eqref{vbb35'}  on
$\O \times (0,T].$

After some minor modifications, one can  check easily that all the estimates in Lemmas \ref{lem10} and \ref{lem11} still hold  for $\al\in [3/4,2)$ and $\ga\in(1,3).$  That is

\begin{lemma}\la{lem3.1}  Let $\al\in [3/4,2)$ and $\ga\in(1,3).$ Then there exists some generic constant $C$  independent of $\ve$ and $T$ such that
\eqref{bb30} and  \eqref{bb5'}   hold with $\O=\mathbb{T}^3.$

\end{lemma}

To obtain the Mellet-Vasseur type estimates for the three-dimensional case, we need to impose  some additional constraints on $\ga$ and $\al.$
\begin{lemma} \la{vlem01} Assume that  $\al\in [3/4,2)$ and that $\ga\in(1,3) $  satisfies $\ga\in ((\al+1)/2,6\al-3) $ in addition.     Then  there exists some generic constant $C$    depending on $T$ but independent of $\ve$ such that \be \la{vbb'01} \sup_{0\le t\le T}\int\n (e+|u|^2)\ln (e+|u|^2)dx\le C.\ee\end{lemma}

{\it Proof.} It is easy to check that \eqref{com1} still holds. Hence, it remains to estimate the righthand side of \eqref{com1}. In fact,
since $\ga\in ((\al+1)/2,6\al-3),$ we have  $ \de=\frac{6\al-3-\ga}{5\ga+3\al-6}\in (0,1).$ Then
 \be \la{wb70} \ba &\ioo  \ln^2 (e+|u|^2)  \n^{2\ga-\al}  dxdt\\ &\le  C \int_0^T\left(\int  \n\ln^{2/\de} (e+|u|^2)dx\right)^\de \left(\int \n^{(2\ga-\al-\de)/(1-\de)}  dx\right)^{ 1-\de }dt\\ &\le  C +C\int_0^T\int \n^{\al-1+5\ga/3 }  dx dt \le  C , \ea\ee where in the last inequality one has used \eqref{e9}. This   together with   \eqref{com1} yields \eqref{vbb'01}. The proof of Lemma \ref{vlem01} is completed.

Furthermore, the following   estimates on the $L^\infty(0,T;L^4(\O))$-norm of $\n^{1/4}u$ will be used later.
\begin{lemma}
\la{vlem12} Assume that $\al\in [3/4,2)$ and that $\ga\in (1,3).$  Then  there exists some constant $C(\ve) $  depending on $\ve$  and $T$ such that   \be\la{ve1}\sup_{0\le t\le T}\int\n  |u|^{4}dx+ \int_0^T\int h_\ve(\n)|u|^2 |\na u|^2dxdt\le C(\ve).\ee Moreover, if in addition $ (\al+1)/2\le\ga\le 3\al-1,$ there exists some constant $C_1$ independent of $\ve$ such that \be\la{k1}\sup_{0\le t\le T}\int\n  |u|^{4}dx+ \int_0^T\int h_\ve(\n)|u|^2 |\na u|^2dxdt\le C_1.\ee\end{lemma}

{\it Proof.} First, multiplying \eqref{vba1}$_2$ by
$   |u|^{2}u$ and  integrating by parts give
\be \la{wb73}\ba & \frac14(\int \n |u|^4dx)_t +  \int h_\ve(\n)|u|^2|\na u |^2dx + 2\int h_\ve(\n) |u|^2 |\na|u||^2dx\\&\quad +   \int g_\ve(\n) (\div u)^2|u|^2dx+ 2\int g_\ve(\n) \div u |u| u\cdot\na |u|dx\\ &\quad+ e^{-\ve^{-3}} \int(\n^{ \ve^{-2}}+\n^{-\ve^{-2}})|u|^{4}dx \\ &=  \frac\ve4 \int \n^{1/2}\div(\n^{- 1/2}h_\ve'(\n)\na \n)|u|^{4}dx+ \int P\div(|u|^2u)dx\\ &= -\frac{ \ve }{8}\int \n^{-1}h_\ve'(\n)|\na\n |^2|u|^{4}dx+\ve  \int h_\ve'(\n) |u|^{3}\na\n\cdot\na |u|dx\\ &\quad+ \int P\div(|u|^{ 2}u)dx\\ &\le  4 \ve  \int \n h_\ve'(\n) |u|^{ 2 }|\na  u|^2dx+ \int P\div(|u|^{ 2}u)dx\\ &\le  4(\al+\ga) \ve  \int  h_\ve (\n) |u|^{ 2 }|\na  u|^2dx+ \int P\div(|u|^{ 2}u)dx. \ea\ee

Next,
  Cauchy's inequality  implies that for any $\ti\al\in [ 3/4,16/5],$
\bnn \la{wini2}\ba &-  2  \int \n^{\ti\al}  |u|^{2}|\na |u||^2dx+(1-\ti\al )\int \n^{\ti\al}(\div u)^2|u|^{2}dx\\&+ 2 (1-\ti\al )\int \n^{\ti\al} \div u |u| u\cdot\na |u|dx\\ &\le    \left(1-\ti\al+\frac{(1-\ti\al)^2 }{2}\right)\int \n^{\ti\al}(\div u)^2|u|^{ 2}dx\\ &\le    \frac{9}{32} \int \n^{\ti\al}(\div u)^2|u|^{ 2}dx \\ &\le    \frac{27}{32}\int \n^{\ti\al}|\na u|^2|u|^{ 2}dx   ,\ea\enn
 which, combined with \eqref{ini9}, shows
\bnn \la{wb74}\ba  &  -2\int h_\ve(\n) |u|^{ 2} |\na|u||^2dx-\int g_\ve(\n) (\div u)^2|u|^{ 2}dx\\&-2\int g_\ve(\n) \div u |u| u\cdot\na |u|dx \\&\le  \frac{27}{32} \int h_\ve(\n)|u|^{ 2} |\na u |^2dx.\ea\enn

Substituting this into \eqref{wb73} yields that there exists some constant $C_1$ independent of $\ve$ such that
 \be  \ba\la{wb75}&\sup_{0\le t\le T}\int\n |u|^{4}dx+\int_0^T\int h_\ve(\n)|u|^{ 2} |\na u |^2dxdt \\&\le C_1+ C_1\int_0^T \|\n^{2\ga-\al-1/2}\|_{L^2(\O)}\left(1+\|\n^{1/4} u \|_{L^4(\O)}^4\right)dt, \ea\ee
where one has used the following estimate
\bnn\ba  &\int P\left|\div(|u|^2u)\right|dx \\ &\le \frac{1}{32}\int\n^\al |u|^2|\na u|^2dx+C_1 \int\n^{2\ga-\al}|u|^2dx\\ &\le \frac{1}{32}\int\n^\al |u|^2|\na u|^2dx+C_1 \|\n^{2\ga-\al-1/2}\|_{L^2(\O)}\left(1+\|\n^{1/4} u \|_{L^4(\O)}^4\right).\ea\enn

Then,   if $ \ga\in [(\al+1)/2,  3\al-1],$
it holds that  \be \la{wbv5} \ba \int_0^T\|\n^{2\ga-\al-1/2}\|_{L^2(\O)}dt&\le C_1\int_0^T\left(\|\n\|^{1/2}_{L^1(\O)}+ \|\n^{\ga+2\al-3/2}\|_{L^2(\O)}\right)dt  \le C_1,\ea\ee    where in the second inequality one  has used
\be \la{k2}\ba  \int_0^T  \|\n^{\ga+2\al-3/2}\|_{L^2(\O)} dt &\le  \int_0^T\|\n\|^{\al-1/2}_{L^{6\al-3}(\O)} \|\n\|^{ \ga+\al-1 }_{L^{3(\ga+\al-1)}(\O)}dt\\&\le C_1+C_1  \int_0^T\|\na \n^{ (\ga+\al-1 )/2}\|^2_{L^{2}(\O)}dt\le C_1 \ea\ee due to \eqref{bb30},  \eqref{bb5'}, and the Sobolev inequality.   It follows from \eqref{wb75}, \eqref{wbv5},  and the Gronwall inequality that \eqref{k1} holds.

Finally, it follows from  \eqref{bb5'} that \bnn\ba \sup_{0\le t\le T}\|\n^{2\ga-\al-1/2}\|_{L^2(\O)} &\le C (\ve)  ,\ea\enn which together with \eqref{wb75}  and the Gronwall inequality gives \eqref{ve1}. The proof of Lemma \ref{vlem12} is thus completed.

With \eqref{bb30},  \eqref{bb5'}, and   \eqref{ve1} at hand, similar to   Lemmas \ref{lem13} and \ref{lem13'}, one can prove
\begin{lemma}
\la{vlem13}Assume that $\al\in [3/4,2)$ and that $\ga\in (1,3).$ Then there exists some constant $C $ depending on $\ve$ and $T$ such that for all $(x,t)\in \O\times (0,T)$ \be\nonumber\la{vc1} C^{-1}\le \n(x,t)\le C . \ee
  Moreover, for any $p> 2,$ there exists some positive constant $C $ depending on $\ve,p,$ and $T$ such that \be\nonumber\la{vbc4}   \int_0^T\left(\|(\n,u)_t\|^p_{L^p(\O )}+\|(\n,u)\|_{W^{2,p}(\O)}^p\right)dt\le C . \ee\end{lemma}


\subsection{Compactness results}
Throughout this subsection, it will be always assumed that $\al$ and $\ga$ satisfy the conditions listed in Theorem \ref{vth2}.

We first construct the initial data. Choose $\n_{0,\ve}$   as in \eqref{pd7}.  Hence, \eqref{pd8} and \eqref{pd9} also hold.
If $\al\in [3/4,1],$ define $u_{0\ve}$ as in  \eqref{pd21}.
If $\al\in (1,2),$
set  \be \la{vpd21} u_{0\ve}= \n_{0\ve}^{-1/4}w_{0\ve},\ee where $w_{0\ve} \in C^\infty(\O)$ satisfies \bnn \|w_{0\ve} -m_0/\n_0^{3/4}  \|_{L^4(\O)}\le \ve .\enn
 It is easy to check that  \eqref{pd13} and \eqref{pd22}  are still valid for $u_{0\ve}$ as in  \eqref{pd21}. Moreover, \eqref{pd13}, \eqref{pd22}, and \eqref{pd3q}  hold for    $u_{0\ve}$ as in  \eqref{vpd21}.

We then extend $(\n_{0\ve},u_{0\ve} )$ $\O$-periodically to $\r^3.$ Similar to the two-dimensional case,  using Lemma \ref{vlem13}  and the standard parabolic theory \cite{la1}, one can show that    the problem  \eqref{vba1}  \eqref{ini9}  \eqref{ma1} \eqref{vbb35'}   has  a unique strong solution   $( \n_\ve,u_\ve) $  satisfying for any  $T>0$ and any $p > 2,$\bnn  \nv,\,\,\uv , \,\, \n_{\ve t} ,\,\,u_{\ve t}  ,\,\, \na^2\nv,\,\,\na^2\uv \in L^p(\O\times (0,T)).\enn

Lemma \ref{lem3.1}  thus shows that   there exists some generic positive constant $C$ independent of $\ve  $  and $T  $     such that \eqref{d1} and \eqref{d2} still hold with $\O=\mathbb{T}^3.$ Hence,    the combination of  \eqref{d1} with \eqref{d2} implies that   \eqref{e19'} and \eqref{e19} are still valid for  $\O=\mathbb{T}^3.$

Moreover, for $\ve\to 0^+,$ it is easy to check that Lemma \ref{lema2} holds also  for the case that $\O=\mathbb{T}^3.$

The following lemma   deals with the compactness of the momentum.

\begin{lemma}\la{vlem03} If $\al \in [3/4,1],$   there exists a function $m(x,t)\in L^2(0,T;L^{3/2}(\O))$ such that up to a subsequence,
\be \la{ve10}\nv\uv\rt m \mbox{ in }L^2(0,T;L^p (\O)) ,\ee   for all $p\in [1,3/2).$
Moreover, \be \la{ve12}\nv\uv\rt \n u \mbox{ almost everywhere }(x,t)\in\O\times (0,T) ,\ee  where \bnn  u(x,t)\triangleq\begin{cases} m(x,t)/\n (x,t)& \mbox{ for } \n(x,t)>0,\\0,&\mbox{ for } \n(x,t)=0.\end{cases}\enn
 \end{lemma}
{\it Proof.}
 First, since $\al\in [3/4,1],$  it follows from \eqref{d1},   \eqref{d2}, and the Sobolev inequality that
\be\la{ve11}\ba  & \int_0^T\|\na (\nv \uv)\|_{L^1(\O) }^2dt \\ &\le  C\int_0^T\left( \|\nv \|_{L^{2-\al}(\O) }^{2-\al} \|\nv^{\al/2}\na \uv\|_{L^2 (\O)}^2 + \| |\uv| |\na \nv |\|_{L^1 (\O)}^2\right) dt\\&\le  C\int_0^T\left(  \|\nv \|_{L^1(\O) } +\|\na\nv^{\al-1/2} \|_{L^2(\O) }^{2/(2\al-1)} \right)^{2-\al}\|\nv^{\al/2}\na \uv\|_{L^2 (\O)}^2 dt\\& \quad +   C\int_0^T\| |\uv| |\na \nv |\|_{L^1(\O) }^2dt\\&  \le C +   C\int_0^T\| |\uv| |\na \nv |\|_{L^1 (\O)}^2dt  .\ea\ee

Thus, if $ \ga\ge 3\al-1,$ then $1\in [2\al-1,\al+\ga-1] $ due to $\al\in[ 3/4 ,1].$ Hence,   it follows from \eqref{d1} and   \eqref{d2} that \be \la{k3}\ba& \int_0^T\| |\uv| |\na \nv |\|_{L^1 (\O)}^2dt\\&\le \int_0^T\| \nv^{1/2}\uv\|_{L^2(\O)}^2\left(\|\na \nv^{\al-1/2}  \|_{L^2 (\O)}^2+\|\na \nv^{(\ga+\al-1)/2} \|_{L^2 (\O)}^2\right)dt \le C.\ea\ee

For $1< \ga\le 3\al-1,$     \eqref{d1},     \eqref{d2}, \eqref{k1}, and  \eqref{k2} imply that \be\la{k4}\ba  \int_0^T\| |\uv| |\na \nv |\|_{L^1(\O) }^2dt &\le \int_0^T\| \nv^{-\al+5/4}\|^2_{L^4(\O) }\|\nv^{1/4}\uv\|_{L^4(\O) }^2 \|\na \nv^{\al-1/2}  \|_{L^2(\O) }^2 dt \\&\le C\int_0^T \left(\|\n\|^{1/2}_{L^1(\O)}+\|\n^{\ga+2\al-3/2}\|_{L^2(\O)}\right)dt\le C,\ea\ee
 where in the second inequality one has used   $\ga+4\al\ge 4$ due to $\al\in [3/4,1]$ and $\ga>1.$ Putting \eqref{k3} and  \eqref{k4} into \eqref{ve11} leads to \be\la{k5}  \int_0^T\|\na (\nv \uv)\|_{L^1 (\O)}^2dt  \le  C.\ee

Next, similar to \eqref{e'14}, it holds that
\be\la{ve14}\ba &(\nv\uv)_t+\div(\nv\uv\otimes\uv)-\div( \nv^\al  \na\uv)-(\al-1)\na( \nv^\al \div\uv) +\na \nv^\ga\\&=\ve\div(h'_\ve(\nv)\na\nv\otimes \uv) -\frac{\ve}{2}\nv^{-1}h'_\ve(\nv)|\na\nv|^2\uv -\ve h'_\ve(\nv)\na\nv \cdot \na\uv   \\&\quad -e^{-\ve^{-3}}(\nv^{ \ve^{-2}}+\nv^{-\ve^{-2}})\uv+\ve^{1/3}\div((\nv^{7/8}+\nv^{\ti\ga})\na\uv) \\&\quad -\frac{\ve^{1/3}}{8}\na(\nv^{7/8}\div\uv) +\ve^{1/3}(\ti\ga-1)\na(\nv^{\ti\ga}\div\uv).\ea\ee

Then, each term on the righthand side of \eqref{ve14} can be estimated similarly as those of   \eqref{e152}--\eqref{e52} and \eqref{e53}. Moreover, for the terms on the left hand side of \eqref{ve14}, we have
\bnn\la{ve15}\ba & \ioo  \nv|\uv|^2dxdt  +\ioo  \nv^\ga dxdt\le C,\ea\enn
\be \nonumber\la{ve18}\ba   \ioo\nv^\al |\na\uv|dxdt&\le C  \ioo h_\ve(\nv)|\na\uv|^2 dxdt+C\ioo h_\ve(\nv) dxdt  \\&\le C .\ea\ee Hence, \be \la{nk5}\|(\nv\uv)_t\|_{ L^1(0,T;W^{-1,1}(\O) )}\le C.\ee  With \eqref{k5} and  \eqref{nk5} at hand, one can finish the proof of Lemma \ref{vlem03} similarly as that of Lemma \ref{lem02}.

When $\al\in (1,2),$  the following compactness result of $\nv^{(\ga+1)/2}\uv $ is needed.
\begin{lemma} \la{vlem02} Assume that $\al\in(1,2) $ and that \be\la{k6} \ga \in [2\al-1,3\al-1] .\ee Then there exists a function $m(x,t)\in L^2(\O\times(0,T ))$ such that up to a subsequence,
\be \la{ve25}\nv^{(\ga+1)/2}\uv\rt m\mbox{ in }L^2(0,T;L^p (\O)) ,\ee      for all $p\in [1,2).$
Moreover, \be \la{ve25'}\nv^{(\ga+1)/2}\uv\rt \n^{(\ga+1)/2}u\mbox{ almost everywhere }(x,t)\in\O\times (0,T)  ,\ee  where \bnn u(x,t)\triangleq\begin{cases} m(x,t)/\n^{(\ga+1)/2} (x,t)& \mbox{ for } \n(x,t)>0,\\0,&\mbox{ for } \n(x,t)=0.\end{cases}\enn
 \end{lemma}

  {\it Proof.} First, it follows from
     \eqref{d1},   \eqref{d2}, and   \eqref{k1}  that
\be\la{ve27}\ba  &\int_0^T\|\na (\nv^{(1+\ga)/2} \uv)\|_{L^{6/5}(\O) }^2dt \\ &\le  C\int_0^T\|\nv^{(1+\ga-\al)/2}\|_{L^3 (\O)}^2 \|\nv^{\al/2}\na \uv\|_{L^2 (\O)}^2dt\\&\quad+ C\int_0^T\|\nv^{1/4}\|_{L^{12}(\O) }^2\|\nv^{1/4}\uv\|_{L^4(\O) }^2   \|\na \nv^{ \ga  /2}\|_{L^2 (\O)}^2   dt \\&\le C  , \ea\ee
where in the second inequality one has used the Sobolev  inequality, $\al\in (1,2),$ and $\ga \in [2\al-1,3\al-1].$
Thus, the combination of \eqref{ve27}, \eqref{e30}, and the Sobolev inequality shows
\be\la{ve29}\int_0^T\|\na (\nv^{(1+\ga)/2} \uv)\|_{L^{6/5} (\O)}^2dt+\int_0^T\|  \nv^{(1+\ga)/2} \uv \|_{L^2(\O) }^2dt\le C.\ee

Next, note that \eqref{e31} and  \eqref{e32} both still hold. Moreover,  it follows from \eqref{vba1} that
\be\la{ve33} \ba  \nv^{(\ga+1)/2}   (\uv)_t =&- \nv^{(\ga+1)/2}\uv\cdot\na \uv +\div(\nv^{(\ga-1)/2}h_\ve(\nv)\na\uv)\\&-\frac{\ga-1}{2}\nv^{(\ga-3)/2} h_\ve(\nv)\na\nv\cdot\na\uv \\&+\na(\nv^{(\ga-1)/2}g_\ve(\nv)\div\uv)-\frac{\ga-1}{2}\nv^{(\ga-3)/2} g_\ve(\nv)\na\nv\div\uv  \\&-\nv^{(\ga-1)/2}\na\nv^\ga-e^{-\ve^{-3}} (\nv^{ \ve^{-2}+(\ga-1)/2}+\nv^{-\ve^{-2}+(\ga-1)/2}) \uv. \ea\ee
 Using \eqref{e31},   \eqref{e32}, and  \eqref{ve33}, one can prove   Lemma \ref{vlem02} in a similar way  as that of Lemma \ref{lem02}.

Next, as a consequence of Lemmas \ref{lema2},  \ref{vlem01}, \ref{vlem03}, and  \ref{vlem02}, similar to Lemma \ref{lema1}, one can obtain
\begin{lemma}\la{vlema1} Assume that $\al $ and  $\ga$ satisfy the conditions listed in Theorem \ref{vth2}. Then up to a subsequence, \be \la{ve42} \sqrt{\nv}\uv\rt \sqrt{\n}u \mbox{ strongly in } L^2 (\O\times (0,T)),\ee with \be \la{ev42} \sqrt{\n}u \in L^\infty(0,T;L^2(\O )).\ee\end{lemma}

Finally, similar to Lemma \ref{plem1}, one can prove
\begin{lemma}\la{plem1'} Assume that $\al $ and  $\ga$ satisfy the conditions listed in Theorem \ref{vth2}. Then up to a subsequence,\be \la{em54} \nv^\al \na\uv\rt \n^\al\na u \mbox{ in }\mathcal{D}',\ee\be \la{em55} \nv^\al (\na\uv)^{\rm tr}\rt \n^\al(\na u)^{\rm tr}\mbox{ in }\mathcal{D}',\ee\be \la{em56} \nv^\al \div\uv\rt \n^\al\div u \mbox{ in }\mathcal{D}'.\ee\end{lemma}
\subsection{Proof of Theorem  \ref{vth2}:  $\O=\mathbb{T}^3$}

 Using \eqref{d1}, \eqref{d2},  \eqref{e19'},   \eqref{e19}, and Lemmas \ref{lema2}, \ref{vlema1}, and  \ref{plem1'}, one can finish the proof of  Theorem  \ref{vth2} where $\O=\mathbb{T}^3$ in a  similar way as  that for the 2-dimensional periodic case  in Theorem  \ref{th2}.

\section{ Proof of Theorem \ref{qvth2}: $\O=\mathbb{T}^3$}

In this section, since the approximate solutions in the proof  of Theorem \ref{vth2}  cannot be applied directly to the case of system \eqref{ii1} \eqref{aa'1} in three dimension space, we will construct a new approximate system which can be applied to obtain the global weak solutions to the three-dimensional system \eqref{ii1} \eqref{aa'1} with $h=\n$ and $g=0.$

\subsection{A priori estimates }

For constants $p_0 $ and $\ve$ satisfying \bnn\la{qq1}p_0=50,\quad 0<\ve\le \ve_1\triangleq\min\{10^{-10},\eta_0\},\enn with $\eta_0$ as in \eqref{pini1},    we  consider  the following approximate system
      \be\la{qba1} \begin{cases}   \n_t+\div(\n u)  =   \ve v\Delta v +\ve v\div( |\na v|^2 \na v)+\ve\n^{-p_0},
     \\   \n u_t+ \n u\cdot \na u   -\div( \n    \mathcal{D} u) +\na P  \\=  \sqrt{\ve}\div(\n \na u)+ \ve v |\na v|^2 \na v\cdot \na u-\ve \n^{-p_0} u- \ve\n   |u|^{3}u    ,  \end{cases}\ee
  where  $v\triangleq\n^{1/2}. $   The initial conditions of the system \eqref{qba1} are imposed as:
\be \la{qbb35} (\n,u)(x,0)= (\n_{0\ve},u_{0\ve}) ,\ee where  smooth  $\O$-periodic functions $ \n_{0\ve}>0$ and $u_{0\ve} $ satisfy
\be \la{qpd9}\ba&\|\n_{0\ve}\|_{L^1\cap L^{\ga}(\O) } +\|\na  \n_{0\ve}^{ 1/2 } \|_{L^2(\O)}+\ve \|\na  \n_{0\ve}^{1/2 } \|^4_{L^4(\O)}+\ve \|  \n_{0\ve}^{-p_0 } \|_{L^1(\O)} \le C, \ea \ee and \be  \la{qpd22}\int_\O \n_{0\ve}   |u_{0\ve}|^{2+\eta_0}dx\le C,\ee for some constant $C$ independent of $\ve.$

 Let $T>0$ be a fixed time and $(\rho ,u)$   be
a smooth solution to \eqref{qba1}  \eqref{qbb35}  on
$\O \times (0,T]. $
Then, we will establish some necessary a priori bounds
for  $(\rho,u)$. The first one is the energy-type inequality.

\begin{lemma} \la{qlem10} There exists some generic constant $C$ independent of $\ve$  such that
\be\la{qbb30}\ba& \sup_{0\le t\le T}\int(\n |u|^2+\n+  \n^\ga+\ve\n^{-p_0}  )dx   +\int_0^T\int    \n   |\mathcal{D} u|^2 dxdt \\&+\ve\int_0^T\int   \left(   |\na v|^4+|\na v|^2|u|^2 + |\na v|^4|u|^2   + \n^{-p_0}|u|^2+\n |u|^{5}\right) dxdt\\& +\ve^2 \int_0^T    \int \n^{-2p_0-1}dx  dt  \le C,\ea\ee where and throughout  this section, for any $f$,
\bnn \int fdx\triangleq\int_\O fdx.\enn
\end{lemma}

{\it Proof.} First, integrating \eqref{qba1}$_1$ over $\O$ yields
\be \la{qbc6} (\int\n dx)_t +\ve \int ( |\nabla v|^2+ |\nabla v|^4)dx  =\ve \int \n^{-p_0}dx .\ee

Next, multiplying $\eqref{qba1}_2$ by $u,$  integrating by parts, and using $\eqref{qba1}_1,$ we have
\be\la{qbb61}\ba& \frac12 (\int\n |u|^2dx)_t   +\int \n  |\mathcal{D}u|^2  dx +\sqrt{\ve}\int \n  |\na u|^2  dx
+\frac{\ve}{2}\int \n^{-p_0} |u|^2dx\\&\quad+\ve\int\n  |u|^{5}dx +\int   u \cdot \na\n^\ga  dx\\&=\frac{\ve}{2}\int  v\Delta v |u|^2dx+\frac{\ve}{2}\int  v\div (|\na v|^2\na v) |u|^2dx +\ve\int v |\na v|^2 \na v\cdot \na u\cdot udx \\&=-\frac{\ve}{2}\int  |\na v|^2 |u|^2dx-\ve \int v\na v\cdot\na u\cdot udx -\frac{\ve}{2}\int   |\na v|^4 |u|^2dx\\&\le -\frac{\ve}{4}\int  |\na v|^2 |u|^2dx+\frac{\sqrt{\ve}}{2} \int \n|\na u|^2dx-\frac{\ve}{2}\int   |\na v|^4 |u|^2dx  .\ea\ee

Then, to estimate the last term on the left hand side  of \eqref{qbb61}, in a similar way as for \eqref{bini5}, one obtains that for $q\not=1,$
\be\la{qbini5}\ba    \int  u\cdot\na  \n^{q }  dx    &= -\frac{q}{q -1}\int\n^{q -1}\div(\n u) dx\\ &= -\frac{q}{q -1}\int\n^{q -1}(-\n_t + \ve v\Delta v+\ve v\div(|\na v|^2\na v)+\ve \n^{-p_0}) dx \\ &= \frac{1}{   q-1 }(\int\n^{q }dx)_t+\frac{q(2q -1)\ve}{ q -1 }\int\n^{q -1} |\na v|^2(1+ |\na v|^2)dx\\&\quad-\frac{q\ve}{q -1}\int\n^{q -1-p_0}dx ,\ea\ee which,  after choosing $q=-p_0,$ implies
\be \la{qjo1}\ba     &\frac{\ve}{ 6(  p_0+1) }(\int\n^{-p_0 }dx)_t+\frac{  p_0(2p_0+1)\ve^2}{ 6(  p_0+1) }\int\n^{-p_0 -1} |\na v|^2(1+ |\na v|^2)dx\\&\quad+\frac{  p_0\ve^2}{6(  p_0+1)}\int\n^{ -1-2p_0}dx \\&= \frac{\ve}{6 } \int  \n^{-p_0 } \div u  dx\\&\le  \frac{p_0\ve^2}{12(p_0+1)}\int\n^{ -1-2p_0}dx+\frac12\int \n (\mathcal{D} u)^2dx .\ea\ee

Finally, adding  \eqref{qbc6}, \eqref{qbb61}, and \eqref{qjo1} together,
 we obtain \eqref{qbb30} after  using \eqref{qbini5},   Gronwall's inequality, and the following simple fact that \bnn \n^{-p_0+\ga-1}\le \n +\n^{-p_0}.\enn
Hence, the proof of Lemma \ref{qlem10} is finished.

In the same spirit of the BD  entropy estimates   due to Bresch-Desjardins \cite{BreschDesjardinsLin05,BreschDesjardins03,
BreschDesjardins03b,BreschDesjardins02}, the following estimates  also hold.

\begin{lemma}
\la{qlem11} There  exists some generic constant $C$ independent of $\ve$ such that
\be\la{qbb5'}\ba& \sup_{0\le t\le T}\int  \left( |\na v|^2+\ve |\na v|^4   \right)dx   + \int_0^T \int \left( \n   |\na u|^2  +   \n^{\ga-2}   |\na \n|^2\right)dxdt  \\&+\ve \int_0^T\int( (\Delta v)^2+|\na v|^2 |\na^2v|^2)dxdt+\ve^2 \int_0^T\int |\na v|^4|\na^2v|^2dxdt\le C.\ea\ee

\end{lemma}

{\it Proof. } First,  set
\bnn  \la{qlalm1} G  \triangleq    \ve v\Delta v +\ve v\div( |\na v|^2 \na v)+\ve\n^{-p_0},\quad
  \vp(\n)\triangleq \log \n
   .\enn
   Following the same procedure leading to \eqref{bb7}, we can get
\be\la{qbb7}\ba& \frac{1+\sqrt{\ve} }{2}(\int\n^{-1} |\na  \n |^2dx)_t + (\int   u\cdot \na  \n  dx )_t+ \int  \n |\na u|^2dx \\&  \quad + \int P'(\n) \n^{-1}|\na \n|^2dx +(1+\sqrt{\ve} )\int  \n^{-1} G\left( \Delta   \n -\frac12  \n^{-1}|\na  \n|^2 \right)dx \\&=- \int  G\div u dx+2  \int  \n \mathcal{D}u: \na  u   dx  + \ve\int       v |\na v|^2 \na v\cdot \na u\cdot\na \log \n  dx\\&\quad- \ve\int  \n^{-p_0} u\cdot\na \log \n  dx- \ve\int \n   |u|^{3}u  \cdot\na \log \n  dx  \triangleq \sum_{i=1}^5 I_i .\ea\ee

Since \bnn \la{qq2}\Delta   \n -\frac12  \n^{-1}|\na  \n|^2 =2v\Delta v,\enn the last term on the left-hand side of \eqref{qbb7} can be rewritten as
\be\la{qq3}\ba &\int\ \n^{-1} G\left( \Delta   \n -\frac12  \n^{-1}|\na  \n|^2 \right)dx\\&=2\ve\int  (\Delta v+\div (|\na v|^2\na v))\Delta v dx+2\ve\int \n^{-p_0-1/2}\Delta v dx \\&=2\ve\int  (\Delta v)^2dx+2\ve\int    |\na v|^2|\na^2v|^2dx+ \ve\int  |\na|\na v|^2|^2dx \\&\quad+2(2p_0+1)\ve\int \n^{-p_0-1 }|\na v|^2 dx ,\ea\ee
where we have used the following simple fact
\bnn\la{qq3'}\ba    \int   \div (|\na v|^2\na v )\Delta v dx  &=  - \int     |\na v|^2\na v\cdot\na\Delta vdx \\&=   \int    |\na v|^2|\na^2v|^2dx+\frac12  \int  |\na|\na v|^2|^2dx  .\ea\enn
Then we will estimate each $I_i (i=1,\cdots, 5)$ on the righthand side of \eqref{qbb7} as follows:
\be \la{qq5}\ba |I_1|& = \left|\int G \div u dx\right|\\& = \ve \left| \int  v(\Delta v+\div(|\na v|^2\na v))\div u dx+\int   \n^{-p_0} \div u dx \right| \\&\le    \frac{\ve^2}{8}\int(  \Delta v+\div(|\na v|^2\na v))^2dx + \ve^2 \int \n^{-2p_0-1} dx +C\int \n (\div u)^2dx,\ea\ee

\be \la{qq4}|I_2|\le \frac{1}{4 }\int   \n | \na  u|^2   dx +C\int  \n  |\mathcal{D} u|^2 dx,\ee

\be \ba |I_3|&=  \ve \left|\int v |\na v|^2 \na v\cdot \na u\cdot \na \log \n dx\right|\\ &=2\ve\left| \int   |\na v|^2 \na v\cdot \na u\cdot \na   v dx\right|\\&= 2\ve\left|\int \pa_j(|\na v|^2    \pa_j   v \pa_i v)u_i dx\right|\\&\le \frac{\ve}{8}\int|\na v|^2|\na^2v|^2dx+C\ve \int  |\na v|^4| u|^2dx,\ea\ee
\be\ba |I_4|&=\ve\left|\int  \n^{-p_0} u\cdot\na \log \n  dx\right|\\&= \frac{\ve}{p_0}\left|\int  \n^{-p_0} \div u   dx\right|\\&
\le C\ve^2\int \n^{-2p_0-1}dx+C\int \n (\div u)^2dx,\ea\ee

\be\la{qq14} \ba |I_5|& =  \ve \left|\int \n   |u|^{3}u\cdot   \na \log \n dx\right| \\&\le C\ve \int \n |u|^{5}dx+C\ve\int (\n+ \n^{-p_0}) dx +\ve \int |u|^2|\na v|^4dx.\ea\ee

Finally, since $v$ satisfies
\be\la{qjo2}\ba  2v_t  -  \ve  \Delta v -\ve  \div( |\na v|^2 \na v)=-  2u\cdot\na v-v\div u+\ve v^{-2p_0-1},\ea\ee
multiplying \eqref{qjo2} by $\ve  ( \Delta v + \div (|\na v|^2\na v)  $  and integrating the resulting equality over $\O$ lead to
\be\la{qq15}\ba &\ve (\int|\na v|^2dx+\frac12\int |\na v|^4dx)_t+\ve^2\int( \Delta v + \div (|\na v|^2\na v) )^2 dx \\& \quad+  (2p_0+1){\ve^2 } \int v^{-2p_0-2 }( |\na v|^2  +  |\na v|^4  ) dx \\& =  {\ve } \int( \Delta v  + \div (|\na v|^2\na v)  ) v\div u   dx+2{\ve } \int( \Delta v  + \div (|\na v|^2\na v)  )  u\cdot\na v dx\\&\le \frac{\ve^2}{8}\int( \Delta v + \div (|\na v|^2\na v) )^2 dx+C\int \n (\div u)^2dx + \frac{\ve}{8}\int |\Delta v|^2dx\\&\quad+ \frac{\ve}{8} \int|\na v|^2|\na^2v|^2dx +C \ve\int |u|^2|\na v|^2(1+|\na v|^2)dx .\ea\ee

 Adding \eqref{qq15} to \eqref{qbb7},   we obtain \eqref{qbb5'} after using Gronwall's inequality, \eqref{qq3}--\eqref{qq14},  \eqref{qbb30},  and the following simple fact:
\bnn\ba&\int(   \div (|\na v|^2\na v) )^2 dx\\&= \int\pa_j(|\na v|^2\pa_i v) \pa_i(|\na v|^2\pa_j v) dx\\&=\int |\na v|^4|\na^2v|^2dx+\int \pa_j|\na v|^2\pa_iv\pa_i|\na v|^2\pa_jvdx\\&\quad  +2\int \pa_j|\na v|^2\pa_iv |\na v|^2\pa_{ij}vdx\\&=\int |\na v|^4|\na^2v|^2dx+\int(\na v\cdot\na|\na v|^2)^2dx+\int |\na v|^2|\na|\na v|^2|^2dx.\ea\enn
The proof of Lemma \ref{qlem11} is thus finished.

With Lemmas  \ref{qlem10} and  \ref{qlem11} at hand, simliar to Lemma \ref{lem01}, we can prove the following   Mellet-Vasseur type estimate (\cite{MelletVasseur05}).
\begin{lemma} \la{qlem01} Assume that $\ga\in(1,3).$ Then  there exists some generic constant $C$      independent of $\ve$ such that \be \la{qbb'01} \sup_{0\le t\le T}\int\n (e+|u|^2)\ln (e+|u|^2)dx\le C.\ee\end{lemma}

{\it Proof.} First, multiplying \eqref{qba1}$_2$ by
$(1+\ln (e+|u|^2))u $ and integrating lead to
\be \la{qbb67} \ba &\frac{1}{2} \frac{d}{dt}\int\n (e+|u|^2)\ln (e+|u|^2)dx+\int\n  \ln (e+|u|^2)  (|\mathcal{D} u|^2+\sqrt{\ve} |\na u|^2) dx
\\&\le C\int \n  |\na u|^2dx  +C\ve\int\n^{-p_0}dx-\int(1+\ln (e+|u|^2))u\cdot\na \n^\ga dx\\&\quad+\frac{\ve}{2} \int (e+|u|^2)\ln (e+|u|^2) v\Delta v  dx,\ea\ee
where we have used the following simple facts that
\bnn\ba &\frac{\ve}{2} \int (e+|u|^2)\ln (e+|u|^2) v\div(|\na v|^2\na v)  dx\\&\quad+\ve\int v|\na v|^2(1+\ln (e+|u|^2))\na v\cdot\na u\cdot u dx \\&=-\frac{\ve}{2} \int (e+|u|^2)\ln (e+|u|^2)  |\na v|^4 dx\le 0 ,\ea\enn and that \bnn \ba& \frac{\ve}{2}\int \n^{-p_0}(e+|u|^2)\ln (e+|u|^2  ) dx-\ve \int \n^{-p_0}(1+\ln (e+|u|^2) )|u|^2dx\\&\le  C\ve \int\n^{-p_0}dx.\ea\enn

Next, similar to \eqref{bb1x} and \eqref{bb2x}, we obtain that
\be \la{qbb1x} \ba &\left|\int(1+\ln (e+|u|^2))u\cdot\na \n^\ga dx\right|\\&\le  C \int  \ln^2 (e+|u|^2)  \n^{2\ga-1}  dx +C\int \n|\na u|^2 dx , \ea\ee and that
\be  \la{qbb2x} \ba &  \frac{\ve}{2} \int (e+|u|^2)\ln (e+|u|^2)v\Delta v dx \\&\le - \frac{\ve}{8}  \int    |\na v |^2(e+|u|^2)\ln (e+|u|^2)dx+\ve \int\n    |\na u|^2dx\\& \quad+C {\ve}  \int |\na v |^2|u|^2dx+\frac{\sqrt{\ve}}{2}  \int  \n  \ln (e+|u|^2)  |\na u|^2dx . \ea\ee

Finally, it follows from   \eqref{qbb67}--\eqref{qbb2x},   \eqref{qbb30}, \eqref{qbb5'},      and  \eqref{qpd22}  that
\be \la{qcom1}\sup_{0\le t\le T}\int\n (e+|u|^2)\ln (e+|u|^2)dx\le C +\int_0^T \int  \ln^2 (e+|u|^2)  \n^{2\ga-1}  dx dt.\ee
Putting \eqref{wb70} where $\al=1$ into \eqref{qcom1} yields \eqref{qbb'01}. The proof of Lemma \ref{qlem01} is completed.

  Next, we will use  a De Giorgi-type procedure to obtain the following estimates on the lower and upper  bounds  of the density which are the keys  to obtain the global existence of strong solutions to the problem \eqref{qba1}   \eqref{qbb35}.

\begin{lemma}
\la{qlem13}There exists some positive constant $C $ depending on $\ve$ such that for all $(x,t)\in \O\times (0,T)$ \be \la{qc1} C^{-1}\le \n(x,t)\le C . \ee \end{lemma}
{\it Proof.} First, it follows from \eqref{qbb5'}, \eqref{qbb30}, and the Sobolev inequality   that
\be\la{qq6}\ba \sup_{0\le t\le T}\|\n\|_{L^\infty}&=\sup_{0\le t\le T}\|  v\|^2_{L^\infty}\\&\le C\sup_{0\le t\le T}\left(\|  v\|_{L^{2 }}+\|\na v\|_{L^4}\right)^2\le \hat C. \ea\ee

Next, we will use a De Giorgi-type procedure to obtain the     lower   bound  of the density. In fact,   since $w\triangleq v^{-1}$ satisfies
\be \la{pbb81} \ba&2w_t+ 2u\cdot \na w- w\div u+\ve w^{2p_0+3}+2\ve  w^{-1}|\na w|^2+2\ve w^{-5}|\na w|^4\\&= \ve\Delta w+\ve \div (w^{-4}|\na w|^2\na w),\ea\ee
multiplying \eqref{pbb81} by $(w-k)_+$ with $ k\ge \|w(\cdot,0)\|_{L^\infty(\O)} =\|\n_0^{-1/2}\|_{L^\infty(\O)} $ yields that
\be\la{pbb82}\ba & \sup_{0\le t\le T}\int  (w-k)_+^2dx +\ve\int_0^T\int(  |\na(w-k)_+ |^2+w^{-4}|\na(w-k)_+ |^4)dxdt \\&\le C\io w|u||\na (w-k)_+|dxdt+C\io (w-k)_+|u||\na w|dxdt\\&\le C\io 1_{\hat A_k} \n^{-4/3}|u|^{ 4/3} dxdt+\frac{\ve}{2} \io w^{-4} |\na(w-k)_+ |^4dxdt, \ea\ee where $\hat A_k\triangleq \{(x,t)\in \O\times (0,T)|w(x,t)>k\}.$
It follows from H\"older's inequality  and \eqref{qbb30}  that
\be\la{pbb83}\ba &\io 1_{\hat A_k} \n^{-4/3}|u|^{ 4/3} dxdt  \\ &\le \left(\io 1_{\hat A_k}\n^{-24/11}  dxdt\right)^{11/15}\left(\io \n |u|^{5} dxdt\right)^{4/15} \\&\le C\left(\io(\n+\n^{-p_0})  dxdt\right)^{1/15}| \hat A_k |^{2/3 }\\&\le C\hat \nu_k^{2/3},\ea\ee where $\hat \nu_k\triangleq  |\hat A_k|.$ Hence,  putting \eqref{pbb83} into \eqref{pbb82}   leads to
\be\la{pbb84}\ba \|(w-k)_+\|_{L^{10/3}(\O\times (0,T))}^2\le C \hat \nu_k^{2/3 }    ,\ea\ee where   we have used the Sobolev inequality \bnn \|(w-k)_+\|_{L^{10/3}(\O\times (0,T))}^2\le C   \sup_{0\le t\le T}\int (w-k)_+^2dx+ C \int_0^T\int  |\na(w-k)_+ |^2dxdt  .\enn
Thus,  \eqref{pbb84} implies that for $h>k,$   \be \la{qq8}\ba \hat \nu_h \le C(h-k)^{-10/3}\hat\nu_k^{10/9 }\ea\ee due to the following simple fact that
\bnn \|(w-k)_+\|_{L^{10/3}(\O\times (0,T))}^2\ge (h-k)^2 |\hat A_h|^{3/5}.\enn It thus follows from \eqref{qq8} and the  De Giorgi-type lemma \cite[Lemma 4.1.1]{wyw} that there exists some positive constant $C\ge \hat C$ such that   \bnn \sup_{(x,t)\in \O\times (0,T)}\n^{-1}(x,t)\le   C,\enn
which together with \eqref{qq6} gives \eqref{qc1} and finishes the proof of Lemma \ref{qlem13}.

We still need the following lemma concerning   the higher order   estimates on $(\n,u)$ which are necessary  to obtain the   global strong solution to the problem  \eqref{qba1}  \eqref{qbb35}.

\begin{lemma}\la{qlem13'} For any $p> 2,$ there exists some constant $C $ depending on $\ve$ and $p $  such that \be\la{qbc4}   \int_0^T\left(\|(\n_t,\na \n_t,u_t\|^p_{L^p(\O) }+\|(\n,\na \n, u)\|_{W^{2,p} (\O) }^p\right)dt\le C . \ee\end{lemma}

  {\it Proof.} First, it follows from \eqref{qc1},   \eqref{qbb30}, and \eqref{qbb5'} that \be \la{qbb89}\sup_{0\le t\le T}  (\|u\|_{L^2  }+\|\na v\|_{L^2\cap L^4 }) +\int_0^T\int \left( |u |^{5}+|\na v|^4|\na^2v|^2 + |\na u |^2\right)dxdt\le C.\ee

Next, it follows from \eqref{qjo2} that  $v=\n^{1/2}$ satisfies
\be\la{qbb87}2 v_t -\ve \div((1+|\na v|^2)\na v) =-\div(u  v+\na w)- \frac{1}{|\O|}\int ( u\cdot \na v-\ve v^{-2p_0-1}) dx,  \ee where for $t>0,$ $w(\cdot,t)$ is the unique solution to the following problem
\be\la{qbb88}\begin{cases} \Delta w=u\cdot \na v-\ve v^{-2p_0-1}- \frac{1}{|\O|}\int ( u\cdot \na v-\ve v^{-2p_0-1}) dx,& x \in \O, \\ \int wdx=0.\end{cases} \ee
Since \eqref{qbb89}     implies  \be \la{qbb1`}\left| \int u\cdot \na v dx \right| \le C\|u\|_{L^2(\O)}\|\na v\|_{L^2(\O)}\le C,\ee  we obtain   that    $\na w$  satisfies for any $p> 2$
\be \la{qbb90} \ba\|\na w\|_{L^p(\O)}&\le C\|\Delta w\|_{L^{3p/(p+3)}}\\&\le C(p)\|u\|_{L^p(\O)}\|\na v\|_{L^3(\O)}+C(p)\\&\le C(p)\|u\|_{L^p(\O)}+C(p),\ea\ee due to \eqref{qbb88}, \eqref{qbb89}, and \eqref{qc1}.

Setting $$\ti v(x,t)\triangleq v(x,t)+\frac{1}{2|\O|}\int_0^t\int ( u\cdot \na v-\ve v^{-2p_0-1}) dxdt,$$ we get from \eqref{qbb87} that \be \la{qbv87} \begin{cases} 2\ti v_t-\ve\div (|\na \ti v|^2\na \ti v)=\div \ti f,\\ \ti v(x,0)=v(x,0), \end{cases}\ee with $\ti f\triangleq\ve \na\ti v-uv-\na w.$

  Thus,     applying the   $L^p$-estimates  \cite[Theorem 1.2]{am1} (see also\cite{bog1,bor1})    to   \eqref{qbv87}  with periodic data   yields   that for any $p\ge 4 $
\be \la{qbb96}  \ba \int_0^T\|\na v \|_{L^{3p} }^{3p}dt &=\int_0^T\|\na \ti v \|_{L^{3p} }^{3p}dt \\ &\le C(p)\left(1+\int_0^T \|\ti f\|_{L^p }^p dt\right)^2\\ &\le C(p) \left(1+ \int_0^T \|u\|_{L^p }^p dt\right)^2+C(p)\left( \int_0^T \|\na \ti v\|_{L^p }^p dt\right)^2 \\ &\le C(p)+C(p) \left(  \int_0^T \|u\|_{L^p }^p dt\right)^2+\frac{1}{2}\int_0^T\|\na   v \|_{L^{3 p} }^{3 p}dt  ,\ea\ee where we have  used \eqref{qbb1`}, \eqref{qbb90},  \eqref{qc1}, and \eqref{qbb89}.

Next, note that \eqref{qba1}$_2$ implies that $u$ satisfies
\be  \la{qbb93}   u_t   -   (\frac{ 1}{2}+\sqrt{\ve})   \Delta u - \frac{ 1}{2}         \nabla \div u=F,\ee where \be \la{qbg3}\ba F&\triangleq -   u\cdot \na u  -\n^{-1}\na P + \ve v^{-1/2} |\na v|^2 \na v\cdot \na u-\ve \n^{-p_0-1} u- \ve   |u|^{3}u   .\ea\ee
Since \bnn \int \Delta u\cdot\na\div udx=\int |\na \div u|^2dx,\enn multiplying  \eqref{qbb93} by $- 2\Delta u $ and integrating the resulting equality over $\O$ lead to\bnn\ba   &(\|\na u\|_{L^2}^2  )_t+\int ((1+2\sqrt{\ve})| \Delta u|^2+|\na\div u|^2)dx\\&\quad+2\ve\int |u|^3(|\na u|^2+3|\na|u||^2)dx\\&=2\int(  u\cdot \na u  +\n^{-1}\na P - \ve v^{-1/2} |\na v|^2 \na v\cdot \na u+\ve \n^{-p_0-1} u)\cdot \Delta udx\\&\le C\|\Delta u\|_{L^2}\left((\|u\|_{L^5}+\||\na v|^3\|_{L^5})\|\na u\|_{L^2}^{2/5}\|\Delta u\|_{L^2}^{3/5}+\|\na v\|_{L^2}+\|u\|_{L^2}\right)\\&\le \frac12\|\Delta u\|_{L^2}^2+C(\|u\|^5_{L^5}+\||\na v|^3\|^5_{L^5})\|\na u\|_{L^2}^2+C,\ea\enn
where in the last inequality we have used \eqref{qbb89}. This together with Gronwall's inequality, \eqref{qbb89}, and \eqref{qbb96} gives
\be \la{qq11}\sup_{0\le t\le T}\|\na u\|^2_{L^2}+\int \|\na^2u\|_{L^2}^2dt\le C.\ee
It thus follows from this and the Sobolev inequality that \bnn \|u\|_{L^{10}(\O\times (0,T))}+\|\na u\|_{L^{10/3}(\O\times (0,T))}\le C,\enn which together with   \eqref{qbb96}--\eqref{qq11} gives \be \la{qq12}\|u_t\|_{L^{2}(\O\times(0,T))}+ \|\na^2 u \|_{L^{2 }(\O\times(0,T))} +\|F\|_{L^{5/2}(\O\times (0,T))} \le C.\ee
Using \eqref{qq12} and applying the standard $L^p$-estimates  to \eqref{qbb93} \eqref{qbg3}  \eqref{qbb35}   with periodic data yield     that for any $p\ge 2$
  \be  \la{qbb99} \ba  \|u_t\|_{L^{p}(\O\times(0,T))}+ \|\na^2 u \|_{L^{p }(\O\times(0,T))} \le C(p)+C(p)\|F\|_{L^{p}(\O\times(0,T))}  .  \ea\ee
In particular, combining   \eqref{qq12} and \eqref{qbb99}  shows
  \bnn  \la{qbb99'} \ba  \|u_t\|_{L^{5/2}(\O\times(0,T))}+ \|\na^2 u \|_{L^{5/2 }(\O\times(0,T))} \le C  .  \ea\enn
This combined with \eqref{qbb89} and the Sobolev inequality (\cite[Chapter II (3.15)]{la1}) yields that for any $q>2$    \bnn \ba \|  u \|_{L^{q}(\O\times(0,T))}+\|\na u \|_{L^{5}(\O\times(0,T))}\le C(q), \ea\enn
which, together with \eqref{qbb96} and \eqref{qbg3},  gives\bnn \|F\|_{L^{9/2}(\O\times (0,T))} \le C.\enn
Combining this with \eqref{qbb99}  yields
  \bnn   \ba  \|u_t\|_{L^{9/2}(\O\times(0,T))}+ \|\na^2 u \|_{L^{9/2 }(\O\times(0,T))} \le C ,  \ea\enn
which together with the Sobolev inequality (\cite[Chapter II (3.15)]{la1}) shows     \bnn \ba \|  u \|_{L^{\infty}(\O\times(0,T))}+\|\na u \|_{L^{45}(\O\times(0,T))}\le C . \ea\enn
Thus, we get    \bnn \la{qq13}\|F\|_{L^{40}(\O\times (0,T))} \le C,\enn which together with \eqref{qbb99} gives
  \bnn    \ba  \|u_t\|_{L^{40}(\O\times(0,T))}+ \|\na^2 u \|_{L^{40 }(\O\times(0,T))} \le C . \ea\enn  The  Sobolev inequality (\cite[Chapter II (3.15)]{la1}) thus implies \bnn\|\na u\|_{L^{\infty}(\O\times(0,T))}\le C.\enn Then, it holds that for any $p>2 , $
 \be \la{qc3} \|u_t\|_{L^p(\O\times(0,T))}+ \|\na^2 u \|_{L^p(\O\times(0,T))} \le C(p).\ee
With \eqref{qc3}   at hand, one can deduce easily from \eqref{qjo2}  and \eqref{qbb35} that for any $p> 2,$
 \bnn \|\n_t \|_{L^p( 0,T,W^{1,p}(\O))}+ \|\na^2 \n \|_{L^p( 0,T,W^{1,p}(\O))} \le C(p),\enn which, together with \eqref{qc3} and \eqref{qbb89},  gives the desired estimate \eqref{qbc4} and finishes the proof of Lemma \ref{lem13'}.

\subsection{Compactness results}

We first construct the  initial data. Let \be\la{qlk1}\si_0\triangleq 10^{-10}.\ee
Choose
\be\nonumber 0\le\ti \n_{0\ve}\in C^\infty (\O),\quad  \|\na\ti\n_{0\ve}^{1/2} \|_{L^4}^4\le \ve^{-4\si_0}   \ee  satisfying \be\nonumber\la{qpd17} \| \ti\n_{0\ve}-\n_0\|_{L^1 (\O)}+ \| \ti\n_{0\ve}-\n_0\|_{L^\ga (\O)}+\|\na ( \ti\n_{0\ve}^{ 1/2}-\n_{0}^{ 1/2})\|_{L^2(\O)}<\ve.\ee Set \bnn\la{qpd7}\n_{0\ve}=\left(\ti\n_{0\ve}^6+\ve^{24\si_0  }\right)^{1/6}.\enn It is easy to check that
\be\la{qpd8} \lim_{\ve\rightarrow 0}\|\n_{0\ve}  -\n_0\|_{L^1(\O)}=0\ee and that there exists some constant $C$ independent of $\ve$ such that \eqref{qpd9} holds.   Define $u_{0\ve}$ as in \eqref{pd21}.
 It is easy to check that  \eqref{pd13} and \eqref{qpd22}  are still valid.

Extend then  $(\n_{0\ve},u_{0\ve} )$ $\O$-periodically to $\r^3.$   The standard   parabolic  theory  \cite{la1}, together with Lemmas \ref{qlem13} and \ref{qlem13'}, thus yields that
the problem \eqref{qba1}  \eqref{qbb35},   where the initial data  $(\n_{0 },u_{0 })$ is replaced by $(\n_{0\ve},u_{0\ve}),$   has  a unique strong solution   $( \n_\ve,u_\ve) $    satisfying \bnn  \nv,\,\uv,\,(\nv)_t ,\,\na(\nv)_t , \,(\uv)_t ,\,\na^2 \nv  ,\,\na^3 \nv  ,\,\na^2 \uv  \in L^p(\b\times(0,T)),\enn for any $T>0$ and  any $p>2.$ Moreover,  all estimates obtained by Lemmas \ref{qlem10}-\ref{qlem01} still hold for $(\nv,\uv).$

Letting $\ve\to 0^+,$ we will modify  the compactness results in Section \ref{qsec1} to prove that   the limit (in some sense) $(\n ,\sqrt{\n }u)$ of $(\nv,\sqrt{\nv}\uv) $ (up to a subsequence) is a weak solution to
\eqref{ii1} \eqref{aa'1} \eqref{hgvv1} \eqref{en1} with $\al=1$ and $\ga\in (1,3).$ We begin with the following strong convergence of $\nv.$
\begin{lemma} \la{qlema2}   There exists a function $\n\in L^\infty(0,T;L^1(\O)\cap L^\ga(\O))$ such that up to a subsequence,
 \be \la{qe3}\nv \ro \n  \mbox{   in }L^\ga (\O\times(0,T )).\ee   \end{lemma}
  {\it Proof.} First, for   $\vf\triangleq \nv^{1/2},$ it follows from \eqref{qbb30} and  \eqref{qbb5'}  that there exists some generic positive constant $C$ independent of $\ve  $      such that
\be\la{qd1}\ba& \sup_{0\le t\le T}\int  (\n_\ve |u_\ve|^2+\n_\ve+\n_\ve^\ga+\ve\nv^{-p_0})dx    +\ioo \n_\ve |\na u_\ve|^2 dxdt \\&+\ve\ioo \left(|\na\vf|^2|\uv|^2+|\na\vf|^4|\uv|^2+\nv^{-p_0}|\uv|^2+\nv|\uv|^5\right) dx  dt \\&+ \ve^{2} \ioo       \nv^{-2p_0-1}dx  dt  \le C ,\ea\ee
and that
\be\la{qd2}\ba& \sup_{0\le t\le T} \int (|\na\vf|^2+\ve|\na\vf|^4) dx   + \ioo  \n_\ve^{ \ga-2}   |\na \n_\ve|^2dxdt \\&+\ve \ioo   \left( (\Delta \vf)^2+ |\na\vf|^2|\na^2\vf|^2 +\ve|\na\vf|^4|\na^2\vf|^2 \right)  dxdt  \le C.\ea\ee

Then,  \eqref{qd1} and \eqref{qd2} yield that
\be \la{qe4}\sup_{0\le t\le T}\|\na \nv \|_{L^{2\ga/(\ga+1)}(\O)  }\le C\sup_{0\le t\le T}\|\nv^{1/2}\|_{L^{2\ga} (\O)  }\sup_{0\le t\le T} \|\na \nv^{ 1/2}\|_{L^2 (\O)  }\le C,\ee and that
\be\la{qq7}\ba \ve^{ 4/3}\int_0^T\|\na \vf\|_{L^6}^6dt&\le \int_0^T\|\na \vf\|_{L^2}^{2/3}(\ve\|\na \vf\|^4_{L^4})^{1/3}(\ve\|\na \vf\|_{L^{12}}^4)dt\\ & \le C\ve\int_0^T\||\na \vf||\na^2\vf|\|_{L^2}^2dt\le C.\ea\ee  Moreover, note that $\nv $ satisfies
\be \la{qe5}\ba (\n_\ve )_t+\div (\nv u_\ve)  =  \ve \vf\Delta \vf+ {\ve} \div(\vf|\na\vf|^2\na \vf) -\ve |\na \vf|^4+\ve\nv^{-p_0} . \ea\ee It follows from \eqref{qd1},     \eqref{qd2}, and    \eqref{qq7}  that
\be \la{qe6}\ba  \sup_{0\le t\le T}\|\n_\ve  u_\ve\|_{L^1 (\O)  }  \le C\sup_{0\le t\le T}\|\nv^{ 1/2}\|_{L^2 (\O)  }\sup_{0\le t\le T}\| \nv^{1/2}\uv\|_{L^2 (\O)  } \le C ,\ea\ee \be \la{qe7}\ba & \ve\ioo (\vf|\Delta \vf|+\vf|\na\vf|^3+|\na \vf|^4)dxdt \\&\le C \ve\int_0^T  (\|\vf\|_{L^2}\|\Delta\vf\|_{L^2}  +\|\vf\|_{L^2}\|\na\vf\|_{L^6}^3+\|\na\vf\|_{L^2}\|\na\vf\|_{L^6}^3) dt \\& \le C\ve^{1/2}+ C \ve^{1/3}\left(\ve^{4/3}\int_0^T\|\na\vf\|_{L^6}^6dt\right)^{1/2} \\& \le C \ve^{1/3},\ea\ee and that \be \la{qe7'}\ba   \ve\ioo \nv^{-p_0}dxdt &  \le    \ve^{1/(2p_0+1)}\left(\ve^{2}\ioo \nv^{-2p_0-1}dxdt\right)^{p_0/(2p_0+1)}\\& \le C \ve^{1/(2p_0+1)}.\ea\ee
The combination of   \eqref{qe5}--\eqref{qe7'} implies that \be\la{qve4}\|(\nv )_t\|_{L^1(0,T;W^{-1,1} (\O) )}\le C.\ee  Since $\ga<3,$ it follows from \eqref{qe4}, \eqref{qve4}, and the
Aubin-Lions lemma  that    \eqref{qe3} holds  for $\ve\to 0^+$ (up to a consequence).
 The proof of Lemma \ref{qlema2} is finished.

Similar to Lemma \ref{vlem03}, we have the following lemma which  deals with the compactness of the momentum.

\begin{lemma}\la{qvlem03} There exists a function $m(x,t)\in L^2(0,T;L^{3/2}(\O))$ such that up to a subsequence,
\be \la{qve10}\nv\uv\rt m \mbox{ in }L^2(0,T;L^p (\O)) ,\ee   for all $p\in [1,3/2).$
Moreover, \be \la{qve12}\nv\uv\rt \n u \mbox{ almost everywhere }(x,t)\in\O\times (0,T) ,\ee  where \bnn u(x,t)\triangleq\begin{cases} m(x,t)/\n (x,t)& \mbox{ for } \n(x,t)>0,\\0,&\mbox{ for } \n(x,t)=0.\end{cases}\enn
 \end{lemma}
{\it Proof.}
 First,    it follows from \eqref{qd1},   \eqref{qd2}, and the Sobolev inequality that
\be\la{qve11}\ba  & \int_0^T\|\na (\nv \uv)\|_{L^1(\O) }^2dt \\ &\le  C\int_0^T\left( \|\nv \|_{L^{1}(\O) }  \|\nv^{1/2}\na \uv\|_{L^2 (\O)}^2 + \| \nv^{1/2} \uv\|_{L^2}^2 \|\na \nv^{1/2}  \|_{L^2 (\O)}^2\right) dt\\&\le  C  .\ea\ee

Next,   it holds that
\be\la{qve14}
 \ba &(\nv\uv)_t+\div(\nv\uv\otimes\uv)-\div(\nv \mathcal{D}\uv) +\na P(\nv)\\&=\ve\vf\Delta\vf \uv+\ve\div( \vf|\na\vf|^2\na\vf\otimes \uv)-\ve |\na\vf|^4\uv \\&\quad +\sqrt{\ve}\div( \nv  \na\uv ) -\ve\nv|\uv|^3\uv.\ea\ee
  For the terms on the left hand side of \eqref{qve14}, we have
\be\la{qve15}\ba & \ioo  \nv|\uv|^2dxdt  +\ioo  \nv^\ga dxdt\le C,\ea\ee
\be \la{qve18}\ba   \ioo\nv  |\na\uv|dxdt&\le C  \ioo \nv |\na\uv|^2 dxdt+C\ioo  \nv  dxdt   \le C .\ea\ee
Moreover, using \eqref{qd1} and \eqref{qd2}, we can estimate  each term on the righthand side of \eqref{qve14}  as follows:
\be\la{qqve1}\ba &\ve\ioo( \vf|\Delta\vf|| \uv| + \vf |\na\vf |^3| \uv|+ |\na\vf|^4|\uv|)dxdt\\&\le C\ve  \int_0^T\|\vf\uv\|_{L^2}\left(\| \Delta \vf\|_{L^2}   +  \|\na \vf \|^3_{L^6} \right)dt\\&\quad+C\left(\ve\ioo  |\na\vf|^4|\uv|^2dxdt \right)^{1/2}\left(\ve\ioo  |\na\vf|^4 dxdt \right)^{1/2}\\&\le C\ve^{1/6},\ea\ee where in the second inequality we have used \eqref{qq7} and \eqref{qe7},
\be\la{qqve2}\ba &\ve\ioo \nv|\uv|^4dxdt\\&\le C\ve^{1/5}\left( \ioo   \nv dxdt \right)^{1/5}\left(\ve\ioo   \nv|\uv|^5dxdt \right)^{4/5}\\&\le C\ve^{1/5}.\ea\ee

 Hence, \be \la{qnk5}\|(\nv\uv)_t\|_{ L^1(0,T;W^{-1,1}(\O) )}\le C.\ee

  With \eqref{qve11} and  \eqref{qnk5} at hand, one can finish the proof of Lemma \ref{qvlem03} similarly as that of Lemma \ref{lem02}.

Next, as a consequence of Lemmas \ref{qlema2}, \ref{qvlem03}, and  \ref{qlem01},  similar to Lemma \ref{lema1}, one can obtain
\begin{lemma}\la{qlema1}   Up to a subsequence, \be \la{qe42} \sqrt{\nv}\uv\rt \sqrt{\n}u   \mbox{ strongly in } L^2 ( 0,T;L^2  (\O)),\ee with \be \la{qe'42} \sqrt{\n}u \in L^\infty(0,T;L^2(\O )).\ee\end{lemma}

Finally, similar to Lemma \ref{plem1}, one can prove  the following convergence of the diffusion terms.
\begin{lemma}\la{qplem1} Up to a subsequence, \be \la{qe54} \nv^\al \na\uv\rt \n^\al\na u \mbox{ in }\mathcal{D}',\ee\be \la{qe55} \nv^\al (\na\uv)^{\rm tr}\rt \n^\al(\na u)^{\rm tr}\mbox{ in }\mathcal{D}',\ee\be \la{qe56} \nv^\al \div\uv\rt \n^\al\div u \mbox{ in }\mathcal{D}'.\ee\end{lemma}
\subsection{Proof  of Theorem  \ref{qvth2}: $\O=\mathbb{T}^3$}
On the one hand,  for any test function $\psi,$   multiplying \eqref{qe5} by $\psi$, integrating the resulting equality over $\O\times (0,T),$ and taking $\ve\rt 0$ (up to a subsequence), one can verify easily after using \eqref{qe3}, \eqref{qe42},  \eqref{qpd8}, \eqref{qe7}, and \eqref{qe7'}  that $(\n,\sqrt{\n}u)$ satisfies \eqref{fin1}.

On the other  hand, let $\phi$ be a test function. Multiplying \eqref{qve14} by $\phi,$   integrating the resulting equality over $\O\times (0,T),$  and taking $\ve\rt 0$ (up to a subsequence),  by Lemmas  \ref{qlema2},  \ref{qlema1},  and \ref{qplem1},  we obtain after using \eqref{qve18}--\eqref{qqve2}     that $(\n,\sqrt{\n}u)$ satisfies \eqref{fin2}. The proof  of Theorem  \ref{qvth2} in the case $\O=\mathbb{T}^3$ is completed.

 \section{Proofs of Theorems \ref{th2},   \ref{vth2},  and \ref{qvth2}:   Cauchy Problem }

 Finally, in this section, we indicate how to generalize the approaches in the previous two sections to deal with the Cauchy problems in the whole spaces.  We start with the 2-dimensional case.

{\it Proof of Theorem  \ref{th2}: $\O=\r^2.$  }  For $\si_0$ as in \eqref{lk1} and $T>0,$  we consider the system \eqref{ba1}--\eqref{ma1}   in $\b\times (0,T)$  with \bnn\la{mac} \b=(-\ve^{-\si_0},\ve^{-\si_0})^2 .\enn

The  initial approximate will be defined as follows. First, choose
\be\nonumber\la{pd5v}\ti \n_{0\ve}\in C^\infty_0 (\b)\cap C^\infty_0 (\O),\quad  0 \le\ti\n_{0\ve} \le \ve^{-4\si_0}   \ee  satisfying \be\nonumber\la{pd17v} \| \ti\n_{0\ve}-\n_0\|_{L^1 (\b)}+ \| \ti\n_{0\ve}-\n_0\|_{L^\ga (\b)}+\|\na ( \ti\n_{0\ve}^{\al-1/2}-\n_{0}^{\al-1/2})\|_{L^2(\b)}<\ve.\ee For $\nu\ge 2$  suitably large such that
$ \nu(\al-1/2)\ge 5,$ set \be\la{pd7v}\n_{0\ve}=\left(\ti\n_{0\ve}^{\nu(\al-1/2)}+\ve^{4\si_0 \nu (\al-1/2)}\right)^{2/(\nu(2\al-1))}.\ee
 It is easy to check  that there exists some positive constant $C$ independent of $\ve$ such that  \eqref{pd9} with $\O$ replaced by $Q_\ve$ still holds. Moreover, it holds that
\be\la{pd8v} \lim_{\ve\rightarrow 0}\|\n_{0\ve}\psi_\ve  -\n_0\|_{L^1(\O)}=0,\ee
where $\psi_\ve \in C_0^\infty(\O)$ satisfies \bnn \psi_\ve(x)=\begin{cases}0,& \mbox{ for }|x|\ge 8\ve^{-\si_0},\\ 1,& \mbox{ for }|x|\le 4\ve^{-\si_0}.\end{cases}\enn

Since $\n_0,m_0$ satisfy \eqref{pini1},   we construct  $w_{0\ve} \in C_0^\infty(\b)\cap C_0^\infty(\O)$ such that \be\nonumber  \|w_{0\ve} - m_0/\n_0^{(1+\eta_0)/ (2+\eta_0)}\|_{L^{2+\eta_0}(\O)}\le \ve .\ee
Set   \be \la{pd21v} u_{0\ve}=  \n_{0\ve}^{-1/(2+\eta_0)}w_{0\ve} .\ee   Then, it holds that
\be\la{pd13v}\lim_{\ve\rightarrow 0}\|\n_{0\ve} u_{0\ve} -m_0\|_{L^1(\O)}=0,\ee  and that \be  \la{pd22v}\int_{Q_\ve} \n_{0\ve}   |u_{0\ve}|^{2+\eta_0}dx\le C.\ee

  Next, let $n=(n_1,n_2)$ denote the unit outward normal to $\pa\b.$ We  impose the initial and boundary conditions on the system  \eqref{ba1}--\eqref{ma1}  as follows:
\be \la{vbb35} \begin{cases}\frac{\pa\n}{\pa n}=0  \mbox{ on } \pa \b, \quad \n (x,0)= \n_{0\ve}(x) ,\quad x\in \b,\\  u\cdot n=0,\pa_1 u_2-\pa_2u_1 =0 \mbox{ on } \pa \b, \quad  u (x,0)=  u_{0\ve}(x) \quad x\in \b.\end{cases}\ee
It follows from  \eqref{vbb35}$_2$ that for any smooth function $f(x) $
\be\la{bv62}\ba   \int_{\pa\b}f\pa_iu_jn_iu_jdS  &=  \int_{\pa\b}f(\pa_iu_j-\pa_ju_i)n_iu_jdS+\int_{\pa\b}f \pa_ju_i n_iu_jdS\\&=-\int_{\pa\b}fu_i u_j\pa_j n_idS=0.\ea\ee

 With \eqref{bv62} at hand, for any $f,$  denoting $$ \int f dx=\int_{\b} f dx,$$ one can check step by step that all the apriori estimates in Lemmas \ref{lem10}--\ref{lem01}   and  \ref{lem13}, where $\O$ is replaced by $Q_\ve,$  still hold for the smooth solution to the problem \eqref{ba1}--\eqref{ma1} \eqref{vbb35}.  It follows from the $L^p$-theory of parabolic system (\cite{dhp}) that Lemma \ref{lem13'}  with $\O$  replaced by $Q_\ve $  also holds. Moreover, for $T>0,$  $p\in (2,\infty),$  and any $F\in L^p(\b\times (0,T)),$ Theorem 2.1 in \cite{dhp} shows that the problem \eqref{bb93}  \eqref{vbb35}$_2$ admits a unique strong solution $u$ on   $\b\times [0,T]$ satisfying \bnn u,u_t,\na^2u\in L^p(\b\times (0,T)),\enn  provided $\n\in C(\ol\b\times [0,T])$ and $\n>0.$ Therefore,   one can use the standard contraction mapping principle    to prove that for any $p> 2 $ and for suitably small $T>0 $ the problem \eqref{ba1}--\eqref{ma1} \eqref{vbb35} has a unique strong solution $(\nv,\uv)$  on $\b\times [0,T]$  satisfying \be \la{qkz1} \nv,\,\uv,\,(\nv)_t ,\,\na(\nv)_t , \,(\uv)_t ,\,\na^2 \nv  ,\,\na^3 \nv  ,\,\na^2 \uv  \in L^p(\b\times(0,T)).\ee Then,  the apriori estimates obtained by Lemmas  \ref{lem13} and \ref{lem13'} yield that the local-in-time strong solution  $(\nv,\uv)$ in fact is a global one, that is, for any $T>0,$ the problem \eqref{ba1}--\eqref{ma1}  \eqref{vbb35} admits a unique strong solution $(\nv,\uv)$ on $\b\times [0,T]$ satisfying \eqref{qkz1} for any $p> 2$. Moreover,  $(\nv,\uv)$  satisfies  all the uniform estimates  (with respect to $\ve$) in Lemmas \ref{lem10}--\ref{lem01}. We then extend $(\nv,\uv)$    to $\O\times[0,T]$ by defining $0$ outside $\b\times [0,T].
 $

 Then after routine modifications of the proofs of Lemmas \ref{lema2}--\ref{lema1}, we conclude  after a  standard diagonal procedure that
\begin{lemma} \la{lema2v}   There exists a function $\n\in L^\infty(0,T;L^1(\O)\cap L^\ga(\O))$ such that up to a subsequence,
\be\la{e2v}  \nv \ro \n  \mbox{   in }L^\ga (0,T; L^\ga_{\rm loc}(\O )).\ee   Moreover, there exists some function $u(x,t)$ such that  \be \la{e'42v} \sqrt{\n}u \in L^\infty(0,T;L^2(\O )),\ee and that up to a subsequence, \be \la{e42v} \sqrt{\nv}\uv\rt \sqrt{\n}u   \mbox{ strongly in } L^2 ( 0,T;L^2_{\rm loc}  (\O)).\ee   \end{lemma}

Finally, it follows from  Lemma \ref{lema2v} and  the proof of   Lemma \ref{plem1}  that    Lemma \ref{plem1} still holds for   $\O=\r^2$.
With Lemmas \ref{lema2v} and \ref{plem1} at hand, after modifying slightly the proof of the periodic case, one can finish the proof of   Theorem  \ref{th2} for the case that $\O=\r^2$. The proof of   Theorem  \ref{th2} is therefore completed.

We now turn to the 3-dimensional case.

{\it Proof of Theorem  \ref{vth2}: $\O=\r^3.$  }  For $\si_0$ as in \eqref{lk1} and $T>0,$    consider the system \eqref{vba1}   \eqref{ini9} \eqref{ma1} in $\b\times (0,T)$  with \be\la{macv} \b=(-\ve^{-\si_0},\ve^{-\si_0})^3 .\ee

Define $\n_{0,\ve}$   as in \eqref{pd7v}.  Hence,  both \eqref{pd9}, where $\O$ is replaced by $Q_\ve,$ and \eqref{pd8v} also hold.
For $\al\in [3/4,1],$ set $u_{0\ve}$ as  in \eqref{pd21v}.
While for $\al\in (1,2),$ let  \be \la{vpd21v} u_{0\ve}= \n_{0\ve}^{-1/4}w_{0\ve},\ee where $w_{0\ve} \in C_0^\infty(\O)\cap  C_0^\infty(\b)$ satisfies \bnn  \|w_{0\ve} - m_0/\n_0^{3/4}  \|_{L^4(\b)}\le \ve .\enn
 It is easy to check that  \eqref{pd13v}  and \eqref{pd22v} are still valid for $u_{0\ve}$ defined in  both cases, \eqref{pd21v} and \eqref{vpd21v}. Moreover, if $\al\in(1,2),$ \be \int_{Q_\ve}\n_{0\ve}|u_{0\ve}|^4dx\le C.\ee

  Next, let $n=(n_1,n_2,n_3)$ denote the unit outward normal to $\pa\b.$ We   impose the initial and boundary conditions on the system  \eqref{vba1}   \eqref{ini9} \eqref{ma1} as follows:
\be \la{vbb35v}\begin{cases}\frac{\pa\n}{\pa n}=0, u\cdot n=0,(\na\times u)  \times n =0 \mbox{ on } \pa \b,  \\ (\n,u)(x,0)=(\n_{0\ve},u_{0\ve}),\quad x\in \b.\end{cases}\ee
Similar to   \eqref{bv62}, by \eqref{vbb35v}$_1,$ it holds that
\be\la{bv62v}   \int_{\pa\b}f\pa_iu_jn_iu_jdS =0,  \ee  for any smooth function $f(x). $

 With \eqref{bv62v} at hand,  denoting that for any $f$,
 $$ \int fdx=\int_{\b} fdx,$$ one can check step by step that all the estimates in Lemmas \ref{lem3.1}--\ref{vlem13},   where $\O$ is replaced by $Q_\ve,$   are still valid for the smooth solution to the problem \eqref{vba1}  \eqref{ini9} \eqref{ma1} \eqref{vbb35v}.  
 Therefore,  similar as that of two-dimensional case,   one can use the standard contraction mapping principle   and the apriori estimates obtained by Lemmas \ref{lem3.1}--\ref{vlem13} to prove that
 the problem \eqref{vba1}  \eqref{ini9} \eqref{ma1} \eqref{vbb35v} has a unique solution $(\nv,\uv)$ on $\b\times [0,T]$  for any $T>0 $ satisfying \eqref{qkz1} for any $p> 2$ and all the uniform  estimates (with respect to $\ve$) in Lemmas \ref{lem3.1}--\ref{vlem12}. We then extend $(\nv,\uv)$    to $\O\times[0,T]$ by defining $0$ outside $\b\times [0,T].
 $

 It then follows from some slight modifications of the proofs of Lemmas \ref{lema2}, and  \ref{vlem03}--\ref{vlema1}, and a  standard diagonal procedure that Lemmas \ref{lema2v} and \ref{plem1'} still hold.
 With the help of these two lemmas,   one can  adapt  the proof of the periodic case to  finish the proof of   Theorem  \ref{vth2} for the case that $\O=\r^3$. The proof of   Theorem  \ref{vth2} is therefore completed.

{\it Proof of Theorem  \ref{qvth2}: $\O=\r^3.$  } First, we choose the  initial approximate  as follows.  For $\si_0$ as in \eqref{qlk1},  let    \be\la{qmacv} \b=(-\ve^{-\si_0},\ve^{-\si_0})^3 .\ee
 Choose
\be\nonumber 0\le\ti \n_{0\ve}\in C^\infty_0 (\b)\cap C^\infty_0 (\O),\quad  \|\na\ti\n_{0\ve}^{1/2} \|_{L^4}^4\le \ve^{-4\si_0}   \ee  satisfying \be\nonumber  \| \ti\n_{0\ve}-\n_0\|_{L^1 (\O)}+ \| \ti\n_{0\ve}-\n_0\|_{L^\ga (\O)}+\|\na ( \ti\n_{0\ve}^{ 1/2}-\n_{0}^{ 1/2})\|_{L^2(\O)}<\ve.\ee Set \be\la{qpdd7}\n_{0\ve}=\left(\ti\n_{0\ve}^6+\ve^{24\si_0  }\right)^{1/6}.\ee It is easy to check that   both \eqref{qpd9}, where $\O$ is replaced by $Q_\ve,$ and \eqref{pd8v}   hold true.
We then choose $u_{0\ve}$ as  in \eqref{pd21v} which satisfies  \eqref{pd13v}  and \eqref{pd22v}.

  Next, let $n=(n_1,n_2,n_3)$ denote the unit outward normal to $\pa\b.$ We   impose the initial and boundary conditions \eqref{vbb35v} on the system  \eqref{qba1}. Note that  \eqref{bv62v} still holds. Moreover,
  since $\na \rho\cdot n=0$ on $\pa \b,$ we have
  \be \la{qvcv}  g(x)\na\n\cdot \na (f(x)\na\n\cdot n)=0,\quad \mbox{ on }\pa\b,\ee for any smooth functions $f(x)$ and $g(x).$

    Denoting that for any $f$,
 $$ \int fdx=\int_{\b} fdx,$$ one can use \eqref{bv62v} and \eqref{qvcv} to check step by step that all the estimates in Lemmas \ref{qlem10}--\ref{qlem13'},   where $\O$ is replaced by $Q_\ve,$   are still valid for the smooth solution to the problem \eqref{qba1} \eqref{vbb35v}.  
 Then,      one can use the standard contraction mapping principle   and the apriori estimates obtained by Lemmas \ref{qlem10}--\ref{qlem13'} to prove that
 the problem \eqref{qba1}  \eqref{vbb35v} has a unique solution $(\nv,\uv)$ on $\b\times [0,T]$  for any $T>0 $ satisfying \eqref{qkz1} for any $p> 2$ and all the uniform  estimates (with respect to $\ve$) in Lemmas \ref{qlem10}--\ref{qlem01}. We then extend $(\nv,\uv)$    to $\O\times[0,T]$ by defining $0$ outside $\b\times [0,T].
 $

 It then follows from some slight modifications of the proofs of Lemmas \ref{qlema2}  and \ref{qlema1}, and a  standard diagonal procedure that Lemmas \ref{lema2v} and \ref{qplem1} still hold.
 With the help of these two lemmas,   one can  adapt  the proof of the periodic case to  finish the proof of   Theorem  \ref{qvth2} for the case that $\O=\r^3$. The proof of   Theorem  \ref{qvth2} is therefore completed.

 \end{document}